\documentclass{aims}
\usepackage{amsmath,amsfonts}
  \usepackage{paralist}
  \usepackage{graphics} %% add this and next lines if pictures should be in esp format
  \usepackage{epsfig} %For pictures: screened artwork should be set up with an 85 or 100 line screen
\usepackage{graphicx}  \usepackage{epstopdf} %This is to transfer .eps figure to .pdf figure; please compile your paper using PDFLaTex or PDFTeXify.
 \usepackage[colorlinks=true]{hyperref}
\hypersetup{urlcolor=blue, citecolor=red}

\def\currentvolume{X}
 \def\currentissue{X}
  \def\currentyear{200X}
   \def\currentmonth{XX}
    \def\ppages{X--XX}

\theoremstyle{definition}

\usepackage[OT2,OT1]{fontenc}
\newcommand\cyr
{
	\renewcommand\rmdefault{wncyr}
	\renewcommand\sfdefault{wncyss}
	\renewcommand\encodingdefault{OT2}
	\normalfont
	\selectfont
}
\DeclareTextFontCommand{\textcyr}{\cyr}

\title[Review on contraction analysis] % Running head is the full title or shortened version of the full title. This will appear at the top of odd pages. Please make sure it fits within the width limit.
      {Review on contraction analysis and computation of contraction metrics} % Only the first word and proper nouns should be capitalized

\author[Peter Giesl, Sigurdur Hafstein and Christoph Kawan]{}


 \email{p.a.giesl@sussex.ac.uk}
 \email{shafstein@hi.is}
 \email{christoph.kawan@gmx.de}

\thanks{$^*$Corresponding author: Peter Giesl}

\begin{document}
	\maketitle

% Enter the first author's name and address:
\centerline{\scshape Peter Giesl$^*$}
\medskip
{\footnotesize
% Enter the address of the first author
 \centerline{Department of Mathematics
 }
   \centerline{University of Sussex}
   \centerline{Falmer, BN1 9QH, UK}
} % Do not forget to end {\footnotesize with the sign }

\medskip

\centerline{\scshape Sigurdur Hafstein}
\medskip
{\footnotesize
 % Enter the address of the second and third authors
 \centerline{The Science Institute, University of Iceland}
   \centerline{Dunhagi 5}
   \centerline{107 Reykjavik, 	Iceland}
}

\medskip

\centerline{\scshape Christoph Kawan}
\medskip
{\footnotesize
	% Enter the address of the second and third authors
	\centerline{Institute of Informatics, LMU Munich}
	\centerline{Oettingenstra\ss e 67}
	\centerline{80538 Munich, Germany}
}

\bigskip

% The name of the associate editor will be entered by AIMS editorial staff.
% "Communicated by the associate editor name" is not needed for special issue.
 \centerline{(Communicated by the associate editor name)}

%The abstract of your paper
%\begin{abstract}
%This is the abstract of your paper and it should not exceed
%\textbf{200} words.
%\end{abstract}
%\documentclass[a4paper,11pt]{article}

% Bilder: bsp_end.mws
%
%\usepackage{color}
%\usepackage{amsmath,amssymb,epsfig}%,showkeys}
%%\usepackage{showlabels}
%\usepackage{url}
%%\usepackage{subfigure}
%%\usepackage{pslatex}
%%\usepackage{apalike}
%%\usepackage{SCITEPRESS}     % Please add other packages that you may need BEFORE the SCITEPRESS.sty package.

%\newtheorem{theorem}{Theorem}[section]
%\newtheorem{definition}[theorem]{Definition}
%\newtheorem{lemma}[theorem]{Lemma}
%\newtheorem{remark}[theorem]{Remark}
%\newtheorem{cond}{Conditions}
%\newtheorem{proposition}[theorem]{Proposition}
%\newtheorem{corollary}[theorem]{Corollary}
%\newtheorem{example}[theorem]{Example}
%\newenvironment{proof}{\par \vspace{0.3cm} \noindent{\sc Proof:} \ignorespaces}%
%{\nolinebreak\hfill $\square$\par \medskip}

% Begin Siggi Macros
\newcommand{\bnull}{{\boldsymbol 0}}
\newcommand{\bff}{{\bf f}}
\newcommand{\bw}{{\bf w}}
\newcommand{\bq}{{\bf q}}
\newcommand{\bc}{{\bf c}}
\newcommand{\bd}{{\bf d}}
\newcommand{\be}{{\bf e}}
\newcommand{\bu}{{\bf u}}
\newcommand{\bh}{{\bf h}}
\newcommand{\ba}{{\bf a}}
\newcommand{\bl}{{\bf l}}
\newcommand{\bp}{{\bf p}}
\newcommand{\bb}{{\bf b}}
\newcommand{\bz}{{\bf z}}
\newcommand{\bZ}{{\bf Z}}
\newcommand{\bx}{{\bf x}}
\newcommand{\bA}{{\bf A}}
\newcommand{\bK}{{\bf K}}
\newcommand{\bB}{{\bf B}}
\newcommand{\bN}{{\bf N}}
\newcommand{\bk}{{\bf k}}
\newcommand{\by}{{\bf y}}
\newcommand{\bv}{{\bf v}}
\newcommand{\bg}{{\bf g}}
\newcommand{\bs}{{\bf s}}
\newcommand{\bt}{{\bf t}}
\newcommand{\bo}{{\bf o}}
\newcommand{\bxi}{{\boldsymbol \xi}}
\newcommand{\bet}{{\boldsymbol \eta}}
\newcommand{\balpha}{{\boldsymbol \alpha}}
\newcommand{\bbeta}{{\boldsymbol \beta}}
\newcommand{\bgamma}{{\boldsymbol \gamma}}
\newcommand{\PS}{\text{\bf PS}}
\newcommand{\ps}{\text{\bf ps}}
\newcommand{\bphi}{{\boldsymbol \phi}}
\newcommand{\bpsi}{{\boldsymbol \psi}}
\newcommand{\eot}{{\begin{flushright} \vspace{-0.7cm} {$ \blacksquare $} \end{flushright}}}
\newcommand{\eod}{{\begin{flushright} \vspace{-0.7cm} {$ \Box $} \end{flushright}}}
\newcommand{\lra}{\longrightarrow}
\newcommand{\ra}{\rightarrow}
\newcommand{\diff}[2]{\frac{d{#1}}{d{#2}}}
\newcommand{\pdiff}[2]{\frac{\partial{#1}}{\partial{#2}}}
\newcommand{\floor}[1]{\lfloor #1 \rfloor}
\newcommand{\ceil}[1]{\lceil #1 \rceil}
\newcommand{\Sym}{\operatorname{Perm}}
\newcommand{\spur}{\operatorname{Spur}}
\newcommand{\Bild}{\operatorname{Bild}}
\newcommand{\con}{\operatorname{con}}
\newcommand{\Ker}{\operatorname{Ker}}
\newcommand{\Dim}{\operatorname{Dim}}
\newcommand{\Dom}{\operatorname{dom}}
\newcommand{\rank}{\operatorname{rank}}
\newcommand{\range}{\operatorname{range}}
\newcommand{\CPWAL}{\operatorname{CPWA}}
\newcommand{\Supp}{\operatorname{supp}}
\newcommand{\scal}[2]{{\langle {#1},{#2} \rangle}}
\newcommand{\inter}{\operatorname{int}}
\newcommand{\Graph}{\operatorname{graph}}
\newcommand{\co}{\operatorname{co}}
\newcommand{\Sign}{\operatorname{sign}}
\newcommand{\Proof}{\noindent {\sc Proof:} \\}
\newcommand{\n}{\nonumber}
\newcommand{\cA}{{\mathcal A}}
\newcommand{\cB}{{\mathcal B}}
\newcommand{\cC}{{\mathcal C}}
\newcommand{\cD}{{\mathcal D}}
\newcommand{\bD}{{\bf D}}
\newcommand{\cE}{{\mathcal E}}
\newcommand{\cF}{{\mathcal F}}
\newcommand{\cG}{{\mathcal G}}
\newcommand{\cH}{{\mathcal H}}
\newcommand{\cI}{{\mathcal I}}
\newcommand{\cJ}{{\mathcal J}}
\newcommand{\cK}{{\mathcal K}}
\newcommand{\cL}{{\mathcal L}}
\newcommand{\cM}{{\mathcal M}}
\newcommand{\cN}{{\mathcal N}}
\newcommand{\cO}{{\mathcal O}}
\newcommand{\cP}{{\mathcal P}}
\newcommand{\cQ}{{\mathcal Q}}
\newcommand{\cR}{{\mathcal R}}
\newcommand{\cS}{{\mathcal S}}
\newcommand{\cT}{{\mathcal T}}
\newcommand{\cU}{{\mathcal U}}
\newcommand{\cV}{{\mathcal V}}
\newcommand{\cW}{{\mathcal W}}
\newcommand{\cX}{{\mathcal X}}
\newcommand{\cY}{{\mathcal Y}}
\newcommand{\cZ}{{\mathcal Z}}

\newcommand{\sA}{{\mathtt A}}
\newcommand{\sB}{{\mathtt B}}
\newcommand{\sC}{{\mathtt C}}

\newcommand{\s}{\big{|}}

\newcommand{\fS}{{\mathfrak S}}
\newcommand{\conv}{\operatorname{conv}}

\newcommand{\dist}{\mathop{\mathrm{dist}}\limits}
\newcommand{\diam}{\mathop{\mathrm{diam}}\limits}
\newcommand{\sign}{\mathop{\mathrm{sign}}\limits}
\newcommand{\supp}{\mathop{\mathrm{supp}}\limits}
\newcommand{\tr}{\mathop{\mathrm{tr}}\limits}
\newcommand{\op}[1]{\mathord{\mathop{\mbox{$#1$}}\limits^\circ}}
\newcommand{\id}{\mathop{\mathrm{id}}\limits}
\newcommand{\graph}{\mathop{\mathrm{graph}}\limits}

% Brackets
\newcommand{\abvallr}[1]{\left|#1\right|}
\newcommand{\abvalb}[1]{\big|#1\big|}
\newcommand{\ang}[1]{\langle #1 \rangle}
\newcommand{\angb}[1]{\big\langle #1 \big\rangle}
\newcommand{\clint}[1]{[#1]}
\newcommand{\clintb}[1]{\big[#1\big]}
\newcommand{\clintlr}[1]{\left[#1\right]}
\newcommand{\opint}[1]{(#1)}
\newcommand{\opintb}[1]{\big(#1\big)}
\newcommand{\opintB}[1]{\Big(#1\Big)}
\newcommand{\opintlr}[1]{\left(#1\right)}
\newcommand{\set}[1]{\{#1\}}
\newcommand{\setb}[1]{\big\{#1\big\}}
\newcommand{\setB}[1]{\Big\{#1\Big\}}
\newcommand{\setlr}[1]{\left\{#1\right\}}
\newcommand{\norm}[1]{\|#1\|}
\newcommand{\normb}[1]{\big\|#1\big\|}
\newcommand{\normlr}[1]{\left\|#1\right\|}

\def\I{ \mathbb {I}}
\def\R{ \mathbb {R}}
\def\C{ \mathbb {C}}
\def\Z{ \mathbb {Z}}
\newcommand{\N}{\mathbb{N}}
\def\cN{\mathcal{N}}

\newcommand{\rmd}{\mathrm{d}}

\newcommand{\vp}{\varphi}

\newcommand{\fa}{\quad \text{for all }\,}

% end Siggi Macros
%\newtheorem{theorem}{Theorem}[section]
%\newtheorem{definition}[theorem]{Definition}
%\newtheorem{lemma}[theorem]{Lemma}
%\usepackage{colortbl}
%\usepackage{fancyhdr}
%\usepackage[top=3.1cm, bottom=3.2cm, left=3.0cm, right=3.0cm]{geometry}
%%\newcommand{\hwnew}[1]{\textcolor{blue}{#1}}
%
%\newcommand{\blue}[1]{\textcolor{blue}{#1}}
%
%\newcommand{\red}[1]{\textcolor{red}{#1}}

%\sloppy
%\allowdisplaybreaks

%\begin{document}

	\begin{abstract}
		Contraction analysis considers the distance between two adjacent trajectories. If this distance is contracting, then trajectories have the same long-term behavior. The main advantage of this analysis is that it is independent of the solutions under consideration. Using an appropriate metric, with respect to which the distance is contracting, one can show convergence to a unique equilibrium or, if attraction only occurs in certain directions, to a periodic orbit.
		
		Contraction analysis was originally  considered for ordinary differential equations, but has been extended to discrete-time systems, control systems, delay equations and many other types of systems. 	
		Moreover, similar techniques can be applied for the estimation of the dimension of attractors and for the estimation of different notions of entropy (including topological entropy).
		
		This review attempts to link the references in both the mathematical and the engineering literature and, furthermore, point out the recent developments and algorithms in the computation of contraction metrics.
	\end{abstract}
	
	\tableofcontents
	
	\section{Introduction}\label{sec:intro}
	
	(Asymptotic) stability is one of the key properties of solutions to ordinary differential equations (ODE). For solutions within the basin of attraction of an attractor, e.g. an asymptotically stable equilibrium, small perturbations have no influence on the long-term outcome. Hence, it is desirable to have sufficient conditions (or certificates) for stability of real-world systems. In control theory, one is interested in designing suitable controllers that ensure stability of a specific solution. One of these sufficient conditions is based on Lyapunov functions, which measure the (decreasing) generalized distance between a point and an attractor such as an equilibrium.
	
	A different way to study stability and attraction, which does not require any knowledge about the position of the attractor, is to measure the evolution of the distance between two trajectories. If this distance decreases, i.e.~contracts, then the long-term behavior of both solutions, and thus of all solutions in a certain set, is the same. This approach is called \emph{contraction analysis}.
	
	There is a large literature on contraction analysis and related topics, and although we have attempted to include many aspects, we will have missed some and apologize. The main purpose of the review is to show links and connections between the different concepts, give an overview over the literature and to present algorithms for the computation of contraction metrics.
	
	The review presents the different aspects of contraction analysis in Section \ref{sec:intro}, starting with an explanation of the concepts in the simplest case and giving a historical overview. In Section \ref{sec:ext}, extensions are considered such as contraction only in certain directions and generalizations to other classes of systems. Finally, Section \ref{sec:num} gives an overview over numerical methods for the computation of contraction metrics.
	
	\subsection{Notation}\label{sec:notation}
	
	We denote by ${\mathcal S}_n$ the set of symmetric matrices in $\mathbb R^{n\times n}$ and by ${\mathcal S}_n^+\subset {\mathcal S}_n$ the subset of positive definite matrices. For matrices $A,B\in {\mathcal S}_n$, we write $A\le B$ if the matrix $A-B$ is negative semidefinite and $A<B$ if $A-B$ is negative definite. We write the identity matrix in any $\R^{n\times n}$ as $I$ and leave it to the reader to work out the dimension from the context. Recall that a matrix $A\in\cS_n$ can be factorized as $A=O^TDO$, where $O\in\R^{n\times n}$ is an orthogonal matrix, i.e.~$O^TO=I$, and $D=\operatorname{diag}(d_1,\ldots,d_n)\in\R^{n\times n}$ is a diagonal matrix. If additionally $A\in\cS_n^+$ we have $d_i> 0$ and  $A^p=O^T\operatorname{diag}(d_1^p,\ldots,d_n^p)O$ for all $p\in\R$.
	
	For a differentiable function $f\colon \mathbb R\times \mathbb R^n\to \mathbb R^{n}$, we denote by $\frac{\partial f}{\partial x}(t,x) \in \mathbb R^{n\times n}$ the Jacobian matrix of partial derivatives with respect to the variable $x\in \mathbb R^n$. The classes $\mathcal K_\infty$ and $\mathcal KL$ of comparison functions are defined as follows: $\alpha\in \cK_\infty$ if and only if $\alpha\colon \R_0^+\to \R_0^+$ is continuous, strictly monotonically increasing, $\alpha(0)=0$, and $\lim_{x\to \infty}\alpha(x)=\infty$; $\beta \in \mathcal KL$ if $\beta \colon \R_0^+\times \R_0^+ \to \R_0^+$ is continuous, $\beta(\cdot,y)\in \cK$ for all $y\in\R_0^+$ and $\beta(x,\cdot)$ is strictly monotonically decreasing with $\lim_{y\to \infty}\beta(x,y)=0$ for all $x\in \R^+$.
	
	Let $f\colon \R\times\R^n \to \R^n$. Then $\phi(t,t_0,x_0)$ denotes the solution $x(t)$ of the initial value problem $\dot{x}=f(t,x)$, $x(t_0)=x_0$; for an autonomous equation, we just write $\phi(t,x_0)$ and assume $t_0=0$.
	For $M\colon \R^n\to \R$ (Lyapunov function) or $M\colon \R^n\to \R^{n\times n}$ (contraction metric), we write $\dot{M}(x)$ for the derivative along solutions of the ODE $\dot{x}=f(x)$, i.e.~its entries $\dot{M}_{ij}(x)$, $i=j=1$
	or $i,j=1,\ldots,n$, are given by
	$$
	\dot{M}_{ij}(x)=\frac{d}{dt}M_{ij}(\phi(t,x))\big|_{t=0}=\nabla M_{ij}(\phi(t,x))\cdot f(\phi(t,x))\big|_{t=0}=\nabla M_{ij}(x)\cdot f(x).
	$$
	Usually, a Lyapunov function is denoted by $V(x):=M_{11}(x)$.
	Similarly, $\dot{M}(t,x)$ for a time-dependent metric $M\colon \R\times \R^n\to \R^{n\times n}$ has the entries
	$$
	\dot{M}_{ij}(t,x)=\frac{\partial M_{ij}}{\partial t}(t,x)+\nabla_x M_{ij}(t,x)\cdot f(t,x).
	$$
	For a vector norm $\|\cdot\|$ on $\R^n$, the induced matrix norm on $\R^{n\times n}$ is defined through $\|A\|:=\sup_{x\neq 0}\|Ax\|/\|x\|=\max_{\|x\|=1}\|Ax\|$. We write $\|x\|_p$ for the standard vector $p$-norms $\displaystyle \|x\|_p:=\left(\sum_{i=1}^n|x_i|^p\right)^{\frac{1}{p}}$ for $1\le p <\infty$ and
	$\displaystyle \|x\|_\infty:=\max_{i=1,\ldots,n}|x_i|$. % Euclidean norm $\|x\|_2 = (x^T x)^{1/2}$.
We will often define a metric by an inner product  $\langle v,w \rangle$ through $\|v\|=\sqrt{\langle v,v\rangle}$ and the induced metric  $\dist(x,y)=\|x-y\|$.

	\subsection{Contraction analysis}
	\label{sec:linear}
	
	A contraction metric is a Riemannian metric, with respect to which the distance between two solutions of a differential equation is decreasing.
	
\begin{figure}[ht]
		\centering
		\includegraphics[width=0.8\textwidth]{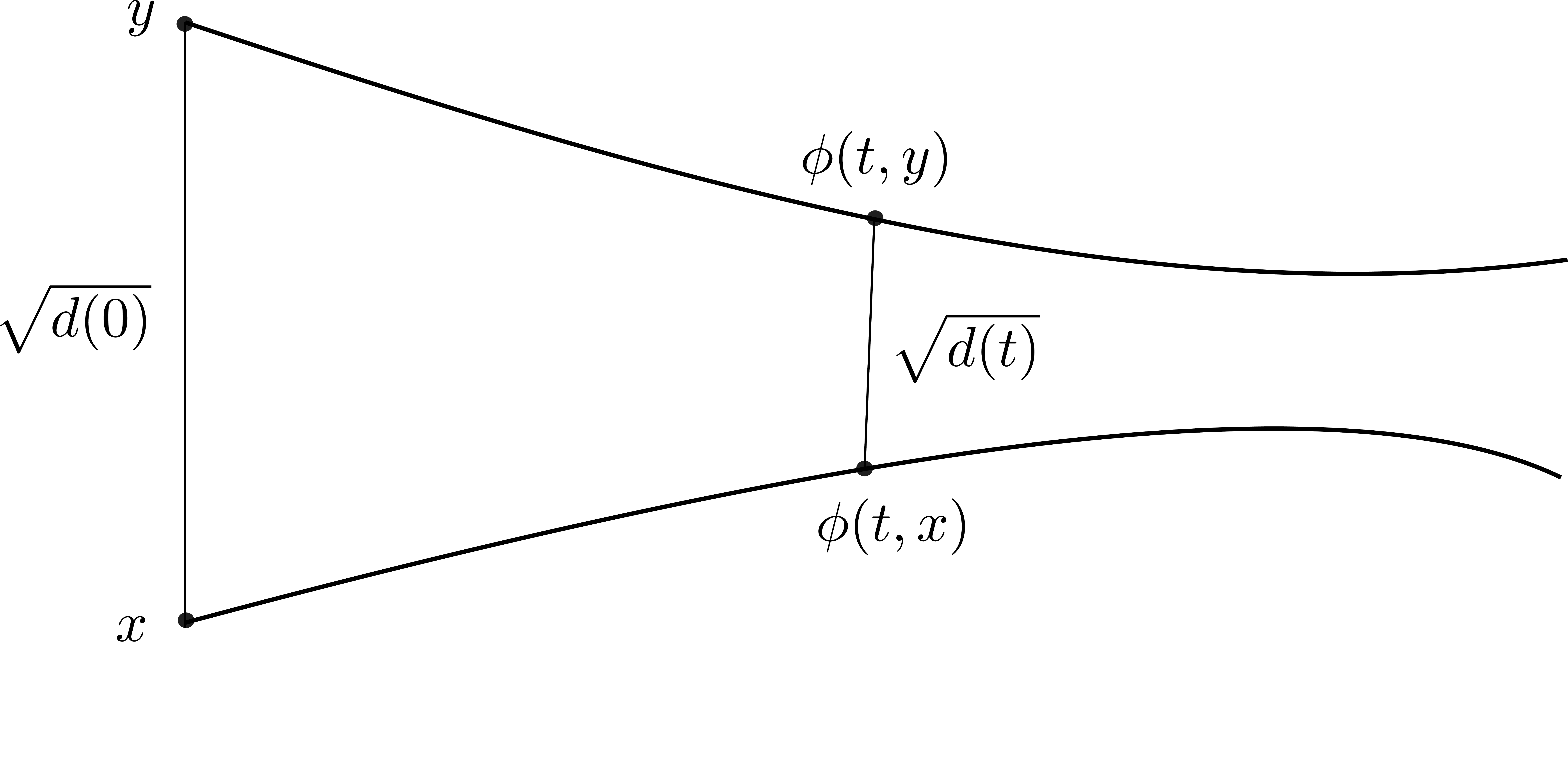}
		\caption{The Euclidian distance $d(t)$ between solutions $\phi(t,x)$ and $\phi(t,y)$ starting at $x$ and $y$ at time $t=0$.}
		\label{fig1}
	\end{figure}

	Let us explain how to derive a condition for the contraction in the simplest case, namely a linear system $\dot{x}=Ax$ with $A\in \mathbb R^{n\times n}$.
	% and denote the solution $x(t)$ of the initial value problem with $x(0)=x_0$ by $\phi(t,x_0)$.
	We consider two solutions $\phi(t,x)$ and $\phi(t,y)$ with initial values $x,y\in\mathbb R^n$ at time $0$, and denote the squared distance between them with respect to the Euclidean metric by%
	\begin{equation*}
		d(t) = (\phi(t,y)-\phi(t,x))^T (\phi(t,y)-\phi(t,x))%
	\end{equation*}
	see Figure \ref{fig1}.
	The derivative of $d(t)$ is given by
	\begin{eqnarray*}
		\frac{d}{dt}d(t)&=& \frac{d}{dt}(\phi(t,y)-\phi(t,x))^T(\phi(t,y)-\phi(t,x))\\
		&&+(\phi(t,y)-\phi(t,x))^T \frac{d}{dt}(\phi(t,y)-\phi(t,x))\\
		&=& (f(\phi(t,y))-f(\phi(t,x)))^T(\phi(t,y)-\phi(t,x))\\
		&&+
		(\phi(t,y)-\phi(t,x))^T (f(\phi(t,y))-f(\phi(t,x)))\\
		%&=&(\phi(t,y)-\phi(t,x))^T A (\phi(t,y)-\phi(t,x))\\
		&=&(\phi(t,y)-\phi(t,x))^T(A^T+ A) (\phi(t,y)-\phi(t,x)).
	\end{eqnarray*}
	If the matrix $A^T+ A\in {\mathcal S}_n$ is negative definite, then $d(t)$ is strictly decreasing (if $x\not= y$) and $\lim_{t\to \infty}d(t)=0$.
	Note that if  $A^T+ A<0$, then there exists $\beta\in \mathbb R^+$ such that
	\begin{eqnarray}
		A^T+ A\le -\beta I\label{ineq0}
	\end{eqnarray}  and we even have
	\begin{eqnarray*}
		\frac{d}{dt}d(t)\le -\beta d(t),
	\end{eqnarray*}
	showing that $d(t)$, and thus $\|\phi(t,y)-\phi(t,x)\|_2$, converges exponentially to zero with rates at least $\beta,\beta/2$, respectively.
	In this case, the time-$t$ map $\phi(t,\cdot)$ for any $t>0$ is contracting and the contraction mapping theorem shows the existence of a unique fixed point, i.e. an equilibrium of the ODE, which is globally attracting.
	Note that the different aspects of this proof will be discussed in detail and more generality later: formula \eqref{ineq0}, the matrix inequality (contraction metric), the decrease of the distance between two solutions (incremental stability) and the convergence of all solutions to a unique solution (convergent dynamics).
	
	These arguments can be generalized by replacing the Euclidean metric $\langle v,w \rangle = v^Tw$  with the constant metric $\langle v,w \rangle_M=v^TMw$, where $M\in{\mathcal S}_n^+$.
	In this case, $$d_M(t)=(\phi(t,y)-\phi(t,x))^TM (\phi(t,y)-\phi(t,x))$$ and the derivative is given by
	\begin{eqnarray*}
		\frac{d}{dt}d_M(t)&=& \frac{d}{dt}(\phi(t,y)-\phi(t,x))^TM(\phi(t,y)-\phi(t,x))\\
		&&		+(\phi(t,y)-\phi(t,x))^TM \frac{d}{dt}(\phi(t,y)-\phi(t,x))\\
		%	&=&(\phi(t,y)-\phi(t,x))^TM (f(\phi(t,y))-f(\phi(t,x)))\\
		%	&=&(\phi(t,y)-\phi(t,x))^T MA (\phi(t,y)-\phi(t,x))\\
		&=&(\phi(t,y)-\phi(t,x))^T(A^TM+M A) (\phi(t,y)-\phi(t,x)).
	\end{eqnarray*}
	Hence, contraction occurs if
	\begin{eqnarray}
		A^TM+ MA<0\label{cont-M0},
	\end{eqnarray} and in this case there exists $\beta\in \mathbb R^+$ with
	\begin{eqnarray}
		A^TM+ MA \le -\beta M,\label{cont-M00}
	\end{eqnarray}  implying even exponential contraction. % {Note that on noncompact sets these conditions are not equivalent.} The subsequent arguments leading to a unique globally stable equilibrium remain the same as above.
	
\begin{figure}[ht]
	\centering
	\includegraphics[width=0.8\textwidth]{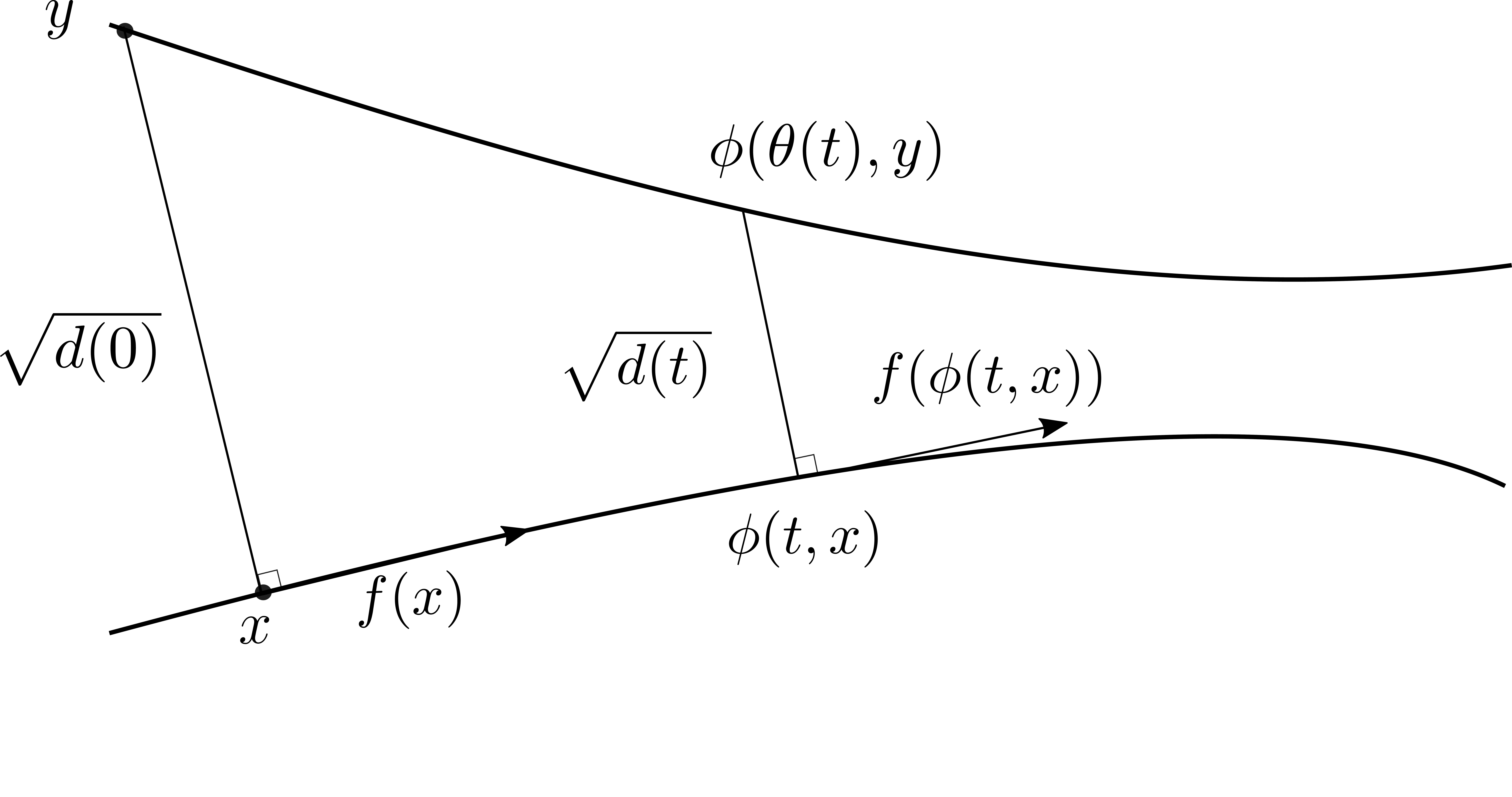}
	\caption{The Euclidian distance $d(t)$ between solutions $\phi(t,x)$ and $\phi(\theta(t),y)$. The time is synchronized with the function $\theta$ such that the difference vector $\phi(\theta(t),y)-\phi(t,x)$ is perpendicular to the flow $f(\phi(t,x))$ at $\phi(t,x)$. }
	\label{fig2}
\end{figure}

	For nonlinear equations, further generalizations are appropriate such as employing a point-dependent metric $\langle v,w\rangle_x=v^TM(x)w$, where $M\in C^1(\mathbb R^n,{\mathcal S}_n^+)$. The (squared) distance $d(t)$ can (i) be measured by the straight line between $\phi(t,x)$ and $\phi(t,y)$ with respect to $M(\phi(t,x))$: $d(t)=(\phi(t,y)-\phi(t,x))^TM(\phi(t,x))(\phi(t,y)-\phi(t,x))$; in this case, a synchronization of the times of the solutions might be applied, i.e.  $$d(t)=(\phi(\theta(t),y)-\phi(t,x))^TM(\phi(t,x))(\phi(\theta(t),y)-\phi(t,x))$$
	with a $\mathcal K_\infty$ function $\theta$.  The function $\theta$ is often chosen such that the difference vector is perpendicular to the flow, either with respect to  the Euclidean metric, i.e.
	$$(\phi(\theta(t),y)-\phi(t,x))^Tf(\phi(t,x))=0,
	$$
	see Figure \ref{fig2},
	or with respect to the metric $M$, i.e., $$(\phi(\theta(t),y)-\phi(t,x))^TM(\phi(t,x))f(\phi(t,x))=0.$$
	This is particularly useful when contraction is only required with respect to certain directions, such as those perpendicular to the flow and all solutions converge to a unique periodic orbit, see Section \ref{sec:period}.

\begin{figure}[ht]
\centering
\includegraphics[width=0.8\textwidth]{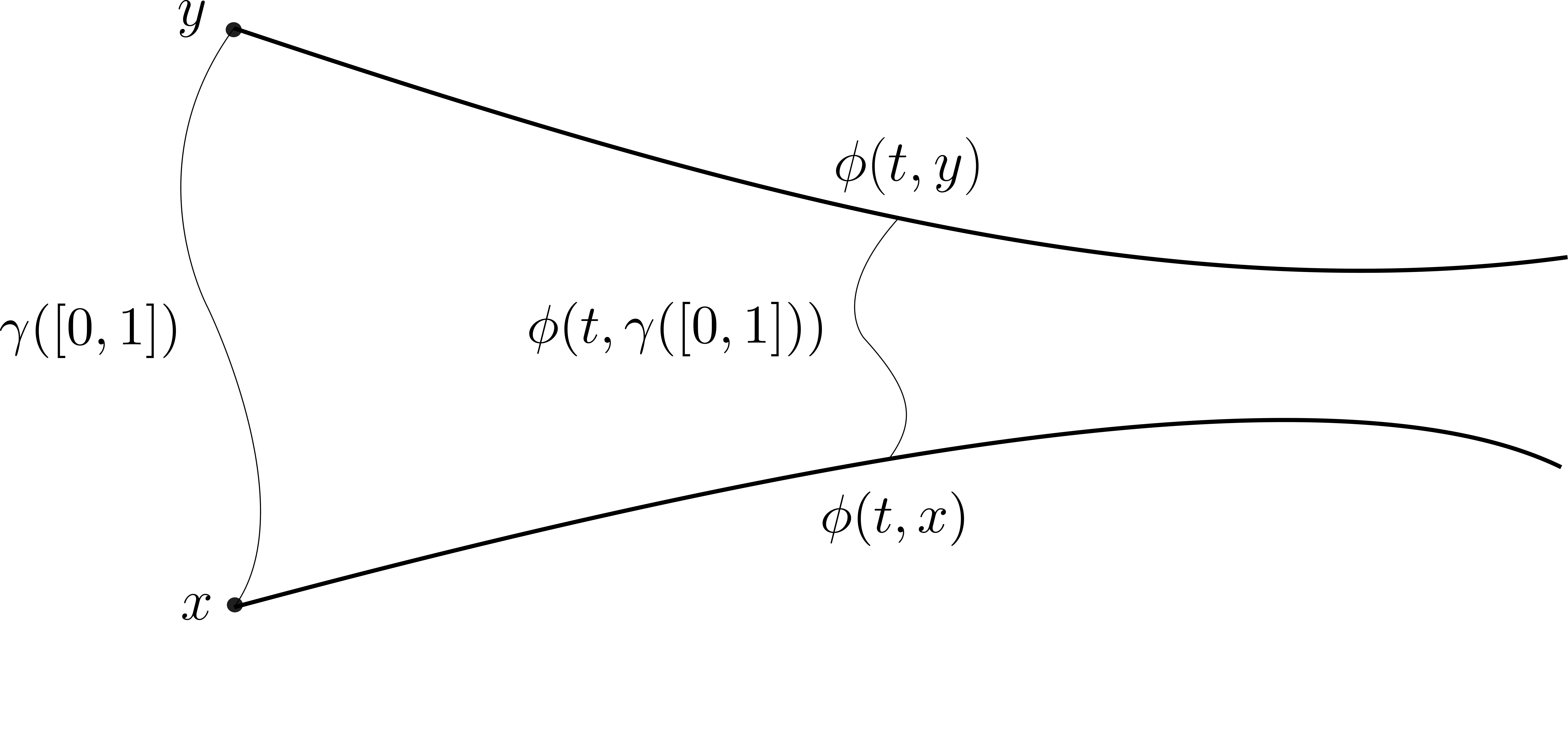}
\caption{The geodesic $\gamma([0,1])$ between $x$ and $y$ and its evolution $\phi(t,\gamma([0,1])$.}
\label{fig3}
\end{figure}

	Another choice (ii) is to consider the distance as the length of the geodesic with respect to the metric $M$; for a schematic depiction see Figure \ref{fig3}. In that case, one  first parameterizes the geodesic between $x$ and $y$ by a $C^1$ function $\gamma\colon [0,1]\to \mathbb R^n$ with $\gamma(0)=x$ and $\gamma(1)=y$. The evolution of this curve under the flow is given by $t\mapsto \phi(t,\gamma([0,1]))$, see Figure \ref{fig3},
	and its length is
	\begin{equation*}
		l(t) = \int_0^1 \Bigl(\frac{\partial}{\partial s}\phi^T(t,\gamma(s))M(\phi(t,\gamma(s)))\frac{\partial}{\partial s}\phi(t,\gamma(s))\Bigr)^{1/2}\,ds.%
	\end{equation*}
	With%
	\begin{equation*}
		l(t,s) := \Bigl(\frac{\partial}{\partial s}\phi^T(t,\gamma(s))M(\phi(t,\gamma(s)))\frac{\partial}{\partial s}\phi(t,\gamma(s))\Bigr)^{1/2},%
	\end{equation*}
	the derivative is given by
	\begin{eqnarray*}
		\frac{d}{dt}l(t) &=& \int_0^1 \frac{1}{2l(t,s)}
		\bigg( \frac{\partial}{\partial s}f^T(\phi(t,\gamma(s)))M(\phi(t,\gamma(s)))\frac{\partial}{\partial s}\phi(t,\gamma(s))\\
		&&+ \frac{\partial}{\partial s}\phi^T(t,\gamma(s))\dot{M}(\phi(t,\gamma(s)))\frac{\partial}{\partial s}\phi(t,\gamma(s))\\
		&&+ \frac{\partial}{\partial s}\phi^T(t,\gamma(s))M(\phi(t,\gamma(s)))\frac{\partial}{\partial s}f(\phi(t,\gamma(s)))\bigg)\,ds\\
		&=&\int_0^1 \frac{1}{2l(t,s)}\bigg( \frac{\partial}{\partial s}\phi^T(t,\gamma(s)) \frac{\partial f^T}{\partial x}(\phi(t,\gamma(s)))M(\phi(t,\gamma(s)))\frac{\partial}{\partial s}\phi(t,\gamma(s))\\
		&& +\frac{\partial}{\partial s}\phi^T(t,\gamma(s))\dot{M}(\phi(t,\gamma(s)))\frac{\partial}{\partial s}\phi(t,\gamma(s))\\
		&&+ \frac{\partial}{\partial s}\phi^T(t,\gamma(s))M(\phi(t,\gamma(s)))\frac{\partial f}{\partial x}(\phi(t,\gamma(s)))\frac{\partial}{\partial s}\phi(t,\gamma(s))\bigg)\,ds\\
		&=& \int_0^1 \frac{1}{2l(t,s)} \frac{\partial}{\partial s}\phi^T(t,\gamma(s))\bigg( \frac{\partial f^T}{\partial x}(\phi(t,\gamma(s)))M(\phi(t,\gamma(s))) + \dot{M}(\phi(t,\gamma(s)))\\
		&&+M(\phi(t,\gamma(s)))\frac{\partial f}{\partial x}(\phi(t,\gamma(s)))\bigg)\frac{\partial}{\partial s}\phi(t,\gamma(s))\,ds.
	\end{eqnarray*}
	Recall from Section \ref{sec:notation} that $\dot{M}(x)\in \mathbb R^{n\times n}$ denotes the derivative along solutions of the ODE, i.e.~the entries of $\dot{M}(x)$ are given by $\dot{M}_{ij}(x)=\nabla M_{ij}(x)\cdot f(x)$.
	%$$\dot{M}_{ij}(x)=\frac{d}{dt}M_{ij}(\phi(t,x))\big|_{t=0}=\nabla M_{ij}(\phi(t,x))\cdot f(\phi(t,x))\big|_{t=0}=\nabla M_{ij}(x)\cdot f(x).$$
	
	If %
	\begin{eqnarray}
		\frac{\partial f^T}{\partial x}(x) M(x)+M(x)\frac{\partial f}{\partial x}(x)+\dot{M}(x) \le -\beta M(x)\label{cont-M}
	\end{eqnarray}
	with $\beta\in\mathbb R^+$, then we can conclude that
	$$	\frac{d}{dt}l(t)\le -\frac{\beta}{2}\, l(t)$$
	and $l(t)\le e^{-\frac{\beta}{2} t}l(0)$. Using that the geodesic distance $d_M(\cdot,\cdot)$ satisfies
	\begin{align*}
		d_M(\phi(t,y),\phi(t,x))&\le l(t)\le e^{-\frac{\beta}{2} t}l(0) = e^{-\frac{\beta}{2}t}d_M(y,x),%
	\end{align*}
	we can conclude exponential contraction with {\it contraction rate} $\beta/2$.%
	
	Note that for the special class of metrics of the form $M(x)=e^{V(x)}I$, condition \eqref{cont-M} becomes
	\begin{eqnarray}
		\frac{\partial f^T}{\partial x}(x)+\frac{\partial f}{\partial x}(x)+\dot{V}(x)I \le -\beta I\label{cont-M-Lyap}
	\end{eqnarray}
	which adds the orbital derivative of the scalar-valued function $V$ to \eqref{cont-M00} for $A=\frac{\partial f}{\partial x}(x)$ and $M=I$, resembling conditions for Lyapunov functions $V$, see also Section \ref{sec:dimension}.
	
	Alternatively, we can (iii) use virtual (infinitesimal) displacements $\delta x$, governed by the first variational equation
	\begin{eqnarray}
		\dot{\delta x} = \frac{\partial f}{\partial x}(x)\delta x,\label{virtual}
	\end{eqnarray}
	which will be explored in the next section.

	\subsection{Different formulations of the contraction property}
	\label{sec:diffformcontr}
	Let us introduce the related concept of a local coordinate change
	$$z=\Theta (t,x)x,$$
	where $\Theta:\R\times\R^n\to \R^{n\times n}$ is a smooth function and  $\Theta(t,x)\in\R^{n\times n}$ is a nonsingular matrix for all $(t,x)$. For generality, we allow for time-dependent coordinate changes. This induces a coordinate change of the displacements
	$$\delta z=\Theta (t,x)\delta x,$$
	leading to the virtual dynamics in the new coordinates below, where $x(t)$ and $\delta x(t)$ are the solution to $\dot{x}=f(x)$ and \eqref{virtual}, respectively,
	\begin{align*}
		\frac{d}{dt} \delta z(t) &= \frac{d}{dt} \Theta (t,x(t))\delta x(t) + \Theta (t,x(t))\frac{d}{dt} \delta x(t)\\
		%&=\delta x(t)\dot{\Theta}(t,x(t))+\Theta (t,x(t))\frac{\partial f}{\partial x}(x(t))\delta x(t)\\
		&= \left(\dot{\Theta}(t,x(t))+\Theta (t,x(t))\frac{\partial f}{\partial x}(x(t))\right)\delta x(t)\\
		&= \left(\dot{\Theta}(t,x(t))+\Theta (t,x(t))\frac{\partial f}{\partial x}(x(t))\right)\Theta^{-1}(t,x(t))\delta z(t),
	\end{align*}
	which can be written compactly as
	\begin{eqnarray}
		\dot{\delta z}&=&F(t,x)\delta z,\text{ where} \label{zformF} \\
		F(t,x)&:=&\left(\dot{\Theta}(t,x)+\Theta(t,x) \frac{\partial f}{\partial x}(x)\right)\Theta^{-1}(t,x). \label{refF}
	\end{eqnarray}
	From
	$$
	\frac{d }{dt}\|\delta z\|_2^2=\dot{\delta z}^T\delta z + \delta z^T\dot{\delta z}=\delta z^T(F^T(t,x)+F(t,x))\delta z
	$$
	it follows that
	exponential contraction is characterized by
	\begin{eqnarray}
		F^T(t,x)+F(t,x) \le -\beta I\label{cont-F}
	\end{eqnarray}
	where $\beta\in\mathbb R^+$.
	To see that \eqref{cont-F} is equivalent to condition \eqref{cont-M}, note that
	for matrices $A,B\in\R^{n\times n}$, $B$ nonsingular, one can easily verify using $y=B^{-1}x$ that $x^TAx\le 0$ for all $x\in\R^n$ is equivalent to  $y^TB^TABy\le 0$ for all $y\in \R^n$.  Thus, we can multiply both sides of \eqref{cont-F} with $\Theta^T(t,x)$ from the left and $\Theta(t,x)$ from the right
	to obtain
	\begin{align*}
		&\dot{\Theta}^T(t,x)\Theta(t,x)+\frac{\partial f^T}{\partial x}(x)\Theta^T(t,x)\Theta(t,x)+\Theta^T(t,x)\dot{\Theta}(t,x)+\Theta^T(t,x)\Theta(t,x)\frac{\partial f}{\partial x}(x)\\
		& \le -\beta \Theta^T(t,x)\Theta(t,x),
	\end{align*}
	which can be written with the time-dependent Riemannian metric $$M(t,x)=\Theta^T(t,x)\Theta(t,x),$$ noting that $\dot{M}(t,x)=\dot{\Theta}^T(t,x)\Theta(t,x)+\Theta^T(t,x)\dot{\Theta}(t,x)$, as
	\begin{eqnarray}
		\frac{\partial f^T}{\partial x}(x)M(t,x) +M(t,x)\frac{\partial f}{\partial x}(x)+\dot{M}(t,x) \le -\beta M(t,x).\label{contraction0}
	\end{eqnarray}
	For a time-independent metric $M(t,x)=M(x)$, this last inequality is \eqref{cont-M}.

	Contractivity can also be characterized by so-called \emph{matrix measures}, also called \emph{logarithmic norms}. They are defined in Lewis 1949 \cite{contr1949lewis}, Dahlquist 1958 \cite{num1958Dahlquist}, Desoer \& Vidyasagar 1975 \cite{cont2009DV} and  Michel 2008 \cite{MiHoLi2008stbdynsyst}.  We follow the presentation in Aminzare \& Sontag 2014 \cite{cont2014AS} and Vidyasagar 2002 \cite{vidya2002}.
	
	Given a vector norm on $\R^n$ with its induced matrix norm, we define the associated matrix measure $\mu$ of a matrix $A\in\R^{n\times n}$ as the directional derivative of the matrix norm at $I$ in direction $A$, i.e.
	$$\mu(A)=\lim_{h\to 0^+}\frac{1}{h}(\|I+hA\|-1).$$
	This limit is known to always exist. For the $\|\cdot\|_2$ vector norm, the induced matrix measure is the maximal eigenvalue of $\frac{1}{2}(A^T+A)$, for the $\|\cdot\|_1$ vector norm, it is given by $\max_{j=1,\ldots,n}\left(a_{jj}+\sum_{i\not=j}|a_{ij}|\right)$, and for  the $\|\cdot\|_\infty$ vector norm, it is $\max_{i=1,\ldots,n}\left(a_{ii}+\sum_{i\not=j}|a_{ij}|\right)$.
	
	For a differentiable $x:\R\to \,\R^n$, we have by Taylor's theorem that
	$$x(t+h)=x(t)+h\dot{x}(t)+o(h),\ \ \ \text{where}\ \ \ \lim_{h\to 0}\frac{\|o(h)\|}{h}=0.$$
	If $x(t)$ is a solution to $\dot{x}=A(t)x$, we have for any norm $\|\cdot\|$ on $\R^n$ that
	$$
	\|x(t+h)\|=\|x(t)+hA(t)x(t)+o(h)\| \le \|I+hA(t)\|\|x(t)\|+\|o(h)\|,
	$$
	which  for $h>0$ implies
	$$
	\frac{\|x(t+h)\|-\|x(t)\|}{h} \le  \frac{\|I+hA(t)\|-1}{h}\|x(t)\|+\frac{\|o(h)\|}{h},
	$$
	i.e.~
	$$
	D^+\|x(t)\|:= \limsup_{h\to 0+} \frac{\|x(t+h)\|-\|x(t)\|}{h} \le \mu(A(t))\|x(t)\|
	$$
	where $D^+$ denotes the upper-right Dini derivative and $\mu$ is the matrix measure with respect to the norm $\|\cdot\|$.
	If $\mu(A(t))\le -b$, this implies $\|x(t)\|\le e^{-bt}\|x(0)\|$, hence, for the solution $\delta z(t)$ to  \eqref{zformF} we obtain
	\begin{equation}
		\label{contraction-log}
		\mu(F(t,x))\le -b \ \ \ \text{for all $t,x$ implies}\ \ \ \|\delta z(t)\|\le e^{-b t}\|\delta z(0)\|,
	\end{equation}
	i.e.~exponential attractivity with \emph{contraction rate} $b$.

	%
	
	%Consider the implication \eqref{contraction-log}.
	If $\mu(F(t,x))\le -b $ holds for all $t\in\R_0^+$ and all $x$ in a convex set $\Omega$,  and $x(t)$ and $y(t)$ are solutions of the ODE that remain in $\Omega$, then
	\begin{equation}
		\label{contrrate}
		\|x(t)-y(t)\|\le e^{-bt}\|x(0)-y(0)\| %\ \ \ \ \text{(recall $b=\beta/2$ contraction rate})
	\end{equation}
	for all $t\ge 0$. Recall that the contraction rate $b$ is equal to $\beta/2$ in the inequalities above.
	The proof considers the ODE for the difference $z(t)=x(t)-y(t)$, namely
	$$\dot{z}(t)=f(x(t))-f(y(t))=\left(\int_0^1 \frac{\partial f}{\partial x}(y(t) + \theta z(t))\,d\theta\right) z(t)$$
	and uses that by the subadditivity of the matrix measure and the convexity of $\Omega$ the matrix measure is bounded by $-b$.
	Concerning the long-term behavior, if $\Omega$ is the entire set and an equilibrium exists, then it must be unique and globally asymptotically stable.
	If $f(t,x)$  is $T$-periodic and $\Omega$ is closed and convex, then there is a unique periodic orbit of period $T$ and every solution in $\Omega$ converges (orbitally) to it.
	
	%It is also possible to use the contraction condition as
	%\begin{eqnarray}
	%\label{contraction-log}\mu\left(\frac{\partial f}{\partial x}(t,x)\right)&\le &-c
	%\end{eqnarray}
	%for all $x\in V$, $t\ge 0$, where $V\subset \mathbb R^n$ and
	It is also possible to bring a constant metric $M$ directly into the contraction condition using an appropriate norm in the definition of $\mu$.
	For $M\in \cS^+_n$, one defines the norm $\|x\|_M:=\sqrt{x^TMx}=\|M^{\frac12}x\|_2$ and the induced  matrix norm with respect to the norm $\|\cdot\|_M$ is
	$$
	\|A\|_M=\sup_{\|x\|_M\neq 0}\frac{\|Ax\|_M}{\|x\|_M}= \sup_{\|x\|_M\neq 0} \frac{\|M^{\frac{1}{2}}Ax\|_2} {\|M^{\frac{1}{2}}x\|_2} =\sup_{\|y\|_2\neq 0}\frac{\|M^{\frac{1}{2}}AM^{-\frac{1}{2}}y\|_2}{\|y\|_2}.
	$$
	It follows that the matrix measure $\mu_M(A)$ using the norm $\|\cdot\|_M$ is the maximal eigenvalue $-b$ of
	$$
	\frac{1}{2}\bigl(M^\frac{1}{2}A^TM^{-\frac{1}{2}}+M^{-\frac{1}{2}}AM^{\frac{1}{2}}\bigr) = \frac{1}{2}M^{-\frac{1}{2}} \left(MA^T+AM\right) M^{-\frac{1}{2}}.
	$$
	In particular, $-b$ is the largest number such that $A^TM+MA \le -2b M$ by the argumentation above  and it follows  that
	%$$
	%\mu_M(A) \le -c\ \ \ \ \Longleftrightarrow  \frac{1}{2}(MA^T+AM) \le -c M.
	%$$
	%In case $\mu=\mu_M$ for a $M\in \cS^+$,
	% from the condition \eqref{contraction-log} with $\Theta(x)\equiv I$ that
	\begin{eqnarray*}
		\mu_M\left(A\right) \le  -b  \ \ \ \Leftrightarrow\ \ \  \frac{1}{2}\left(A^TM+M A\right)\le -b M,
	\end{eqnarray*}
	which is \eqref{cont-M00} with $b=\beta/2$.
	
	%With $M\colon\R^n\to \cS^+_n$ and $\mu_x$ pointwise defined as $\mu_{M(x)}$ one can verify through computations identical to above for $\dot{x}=f(t,x)$ that
	%$$
	%\mu_x\left(\frac{\partial f}{\partial x}(t,x)\right) \le -b\ \ \ \Leftrightarrow\ \ \ \frac{\partial f^T}{\partial x}(t,x)M(x)+M(x) \frac{\partial f}{\partial x}(t,x){\dot M ???}\le -2b M(x).
	%$$
	%
	
	%Finally, all of these ideas can be generalized to non-autonomous systems, where in general the metric will depend on space and time, as well as to infinite-dimensional systems, and also to nonlinear and even non-differentiable functions $f\colon \mathbb R^n\to\mathbb R^n$ (from $f(x)=Ax$) by defining logarithmic Lipschitz constants.

	As a final note on different formulations of the contraction property, note that these ideas can be generalized to non-autonomous systems, where in general the metric will depend on space and time, as well as to infinite-dimensional systems.  Additionally, they can be
	extended to systems with non-differentiable right-hand sides  by defining logarithmic Lipschitz constants.

		For \emph{non-autonomous} systems, the notion of \emph{strictly contracting processes} has been defined in the book
		Kloeden \& Yang 2021 \cite[Definition 7.5]{kloedenYang} as follows: a process $\phi$ on a complete metric space $(X,d_X)$,  e.g.~$\phi(t,t_0,x_0)$ is the solution of the non-autonomous ODE $\dot{x}=f(t,x)$ with initial condition $x(t_0)=x_0$, satisfies a \emph{uniformly strictly contracting property} if there exists $L>0$ such that
		$$d_X(\phi(t,t_0,x_0),\phi(t,t_0,y_0))\le e^{-L(t-t_0)} d_X(x_0,y_0)$$
		holds for all $x_0,y_0\in X$ and $t,t_0\ge 0$. In \cite[Theorem 7.2]{kloedenYang}, it is shown that such a process, under some additional assumptions, has a pullback attractor, which is  also forward attracting, with component sets consisting of single, entire solutions. If $f$ is time-periodic, then the entire solution is also periodic. The main idea of the existence proof is to use the contraction property to show that $\phi(0,t_n,x_n)$ for a decreasing sequence $t_n$ with $\lim_{n\to \infty}t_n=-\infty$ and certain assumptions on $x_n$, is a Cauchy sequence, which thus converges due to the completeness of $X$; this idea goes back to the stochastic case in Caraballo, Kloeden \& Schmalfu\ss \ 2004 \cite{Caraballo2004}.

	\subsection{Contraction analysis vs.~Lyapunov stability theory}
	%Let $f\colon \R_0^+\times\R^n \to \R^n$ and
	Stability measures the distance of solutions to an attractor, which should (exponentially) decrease, while contraction measures the distance of solutions to each other, which should (exponentially) decrease. While the first one requires us to know the attractor, the second one does not. Moreover, contraction is robust under perturbations of the system, even under perturbations of the attractor. Stability is a topological property, independent of the metric, while contractivity is a metric property and depends on the specific metric that is used. In particular, a system can be contracting with respect to one metric, but not with respect to  another.  We first present the simplest case, i.e.~autonomous linear systems, and then we discuss a unified framework for the Lyapunov theory and contraction analysis.
	
	%
	%Stability can be shown through a Lyapunov function. A Lyapunov function for an equilibrium $x_0$ is a scalar-valued function $V\colon N\to\mathbb R$, where $N$ is a neighborhood of $x_0$, that is strictly decreasing along solutions in $N\setminus \{x_0\}$ and satisfies $V(x)>V(x_0)$ for all $x\in N\setminus \{x_0\}$.
	\subsubsection{Autonomous linear systems}\label{sec:aut-lin}
	Let us compare both notions for a linear system $\dot{x}=Ax$. The zero solution is asymptotically stable if and only if $A$ is Hurwitz, i.e.~the eigenvalues of $A$ have negative real part; in this case, it is even exponentially stable and there exists a quadratic Lyapunov function $V(x)=x^TMx$, where the symmetric matrix $M$ is the solution of the Lyapunov equation
	\begin{eqnarray}
		A^TM+MA = -C,\label{Lyap-eq}
	\end{eqnarray}
	for any given symmetric, positive definite matrix $C$. Note that if $A$ is Hurwitz and $C\in {\mathcal S}_n^+$, then this equation has a unique solution $M\in {\mathcal S}_n^+$. This is a converse result, showing that a linear system with an asymptotically stable equilibrium always has a quadratic Lyapunov function.
	
	We will see that such a system is not in general contracting with respect to the Euclidean metric, but it is contracting with respect to a suitable constant Riemannian metric.
	For contraction with respect to the Euclidean metric, we require that $A^T+A$ is negative definite, see Section \ref{sec:linear}. This  implies that $A$ is Hurwitz: to see this, let us assume in contradiction to the statement that $A$ is not Hurwitz, then (real case) there exist $v\in\mathbb R^n\setminus \{0\}$ and $\lambda\ge 0$ with $Av=\lambda v$, which shows $v^T(A^T+A)v=2\lambda v^Tv\ge 0$ in contradiction to $A^T+A$ being negative definite. In the complex case, there exist  vectors $v_1,v_2\in\mathbb R^n$ and $\mu\ge 0$, $\nu\in\mathbb R$ with $A(v_1+iv_2)=(\mu+i\nu)(v_1+iv_2)$ and $v_1+iv_2\not=0$.
	Then %$w=\frac{v_1}{\|v_1\|}+\frac{v_2}{\|v_2\|}\not=0$. Then
	\begin{eqnarray*}
		v_1^T(A^T+A)v_1=2[\mu \|v_1\|^2-\nu v_1^Tv_2]\\
		v_2^T(A^T+A)v_2=2[\mu \|v_2\|^2+\nu v_1^Tv_2]
	\end{eqnarray*}
	and, depending on the sign of $\nu v_1^Tv_2$, at least one of the two expressions is $\ge 0$  in contradiction to $A^T+A$ being negative definite.
	
	But the converse implication is false, as the matrix $A=\left(\begin{array}{rr}-1&a\\0&-1\end{array}\right)$ with $a\ge 2$ shows.  Indeed, the matrix is Hurwitz as the eigenvalue is $-1$, but with $v=(1,1)^T$ we have $v^T(A+A^T)v=
	v^T \left(\begin{array}{rr}-2&a\\a&-2\end{array}\right)v=-4+2a\ge 0$.
	However, if we allow for a more general metric defined by $\langle v,w\rangle_M=v^TMw$, then the system is contracting with respect to this metric if and only if
	$$A^TM+MA<0,$$
	see \eqref{cont-M0};
	this can again be achieved by the solution $M$ of the Lyapunov equation \eqref{Lyap-eq} and is thus equivalent to $A$ being Hurwitz.
	
	Summarizing, asymptotic stability for linear autonomous systems is equivalent to the existence of a Lyapunov function and to the existence of a contraction metric, even a point-independent metric $\langle v,w\rangle_M=v^TMw$, see also  Falsaperla, Giacobbe \& Mulone 2012 \cite{falsaperla2012}. We will later study the question of existence of contraction metrics for more general systems, see Section \ref{sec:converse}.

	\subsubsection{Unified framework}
	We now discuss a unified framework for contraction analysis and Lyapunov stability theory.  For this, it is advantageous to add
	to the non-autonomous system $\dot x=f(t,x)$ with solution $x(t)=\phi(t,t_0,x_0)$ the associated variational equation $\dot{\delta x}=\frac{\partial f}{\partial x}(t,\phi(t,t_0,x_0))\delta x$ with solution $\delta x(t)=\psi(t,t_0,x_0,\delta x_0)$, i.e.~to consider the augmented system
	$$
	\frac{d}{dt}[x,\delta x]= \left[ f(t,x),\frac{\partial f}{\partial x}(t,\phi(t,t_0,x_0))\delta x\right].
	$$
	For simplicity, let us assume that the flow defined by $\dot x=f(t,x)$ evolves on $\cM=\R^n$, i.e.~$f\colon \R\times \R^n\to \R^n$; below, we will discuss ODEs on more general manifolds.

	In both Lyapunov theory and contraction analysis, one studies the derivative of a sufficiently smooth  function  $V\colon \R_0^+\times\R^n\times \R^n \to \R$ along the solution trajectories of this system:
	\begin{align}
		\label{genLya}
		&\dot V(t,x(t),\delta x(t)) := \frac{d}{dt}V(t,x(t),\delta x(t))\\
		&\ \ \ \ =\frac{\partial V}{\partial t}(t,x(t),\delta x(t)) +\frac{\partial V}{\partial x}(t,x(t),\delta x(t))f(t,x(t)) \nonumber\\
		&\ \ \ \ \ \ +\frac{\partial V}{\partial \delta x}(t,x(t),\delta x(t)) \frac{\partial f}{\partial x}(t,x(t))\delta x(t). \nonumber
	\end{align}
	By some abuse of notation, the time-dependence of $x$ and $\delta x$ is usually suppressed and equation \eqref{genLya} is simply written as
	\begin{equation}
		\label{genLya2}
		\dot V(t,x,\delta x) = \frac{\partial V}{\partial t}(t,x,\delta x) +\frac{\partial V}{\partial x}(t,x,\delta x)f(t,x)+\frac{\partial V}{\partial \delta x}(t,x,\delta x) \frac{\partial f}{\partial x}(t,x)\delta x.
	\end{equation}
	%Note that for brevity of the presentation we assume that $V$ is autonomous, i.e.~does not explicitly depend on $t$.
	In the classical Lyapunov theory for the asymptotic stability of, say the zero solution $x(t)\equiv 0$, one considers functions $V$ that do not depend on $\delta x$, i.e.~$V(t,x,\delta x)=V_c(t,x)$, and there is no need for the augmented system.  Classical results are, e.g.~that the global, uniform, asymptotic stability of the zero solution is equivalent to the  existence of a Lyapunov function $V_c(t,x)$ fulfilling $\alpha_1(\|x\|_2) \le V_c(t,x) \le \alpha_2(\|x\|_2)$ and  $\dot{V_c}(t,x)\le -\alpha_3(\|x\|_2)$
	for some $\alpha_1,\alpha_2,\alpha_3\in\cK_\infty$ and all $(t,x)\in \R_0^+\times \R^n$. If $\alpha_i(x)=c_ix^2$  for  constants $c_i>0$, then the equilibrium is even globally uniformly exponentially stable,  cf.~e.g. \cite[Ths.~4.10, 4.14]{Khalil2002}.
	If $f(t,x)$ does not depend on $t$ or is periodic in $t$, then one may assume the same of $V_c(t,x)$. The general principle is to let $V_c(t,x)$ measure the distance to some attractor $A$ and demand that \eqref{genLya2} is strictly negative on $N\setminus A$ for an open neighbourhood $N$ of $A$. Since $t\mapsto V_c(t,x(t))$ is decreasing and $V_c$ obtains its minimum on the attractor, mild assumptions guarantee that all solutions starting in a compact sublevel set $V_c^{-1}((-\infty,c]) \subset N$ converge to the attractor when $t\to \infty$.
	
	In contraction analysis, one typically sets $V(t,x,\delta x) := (\delta x)^T M(t,x) \delta x$, where $M(t,x) \in \cS^+_n$ for all $(t,x)$. Then \eqref{genLya2} can be written as
	\begin{equation}
		\label{genLya3}
		\dot V(t,x,\delta x) = (\delta x)^T \left( \frac{\partial f^T}{\partial x}(t,x) M(t,x) + M(t,x) \frac{\partial f}{\partial x}(t,x)+\dot M(t,x) \right) \delta x.
	\end{equation}
	%where $\dot M(t,x)$ is the symmetric matrix with entries $\left[\dot M(x)\right]_{ij}=\frac{\partial  M_{ij}}{\partial t}(t,x)+\nabla_x M_{ij}(x)\cdot f(t,x)$.
	
	In what follows, we consider again general functions $V(t,x,\delta x)$,  and we let $\alpha_i\in \cK_\infty$ for all $i$.
	Somewhat simplified, one can say that the classical Lyapunov stability theory and contraction analysis both essentially boil down to the decrease condition
	\begin{equation}
		\label{genDEC}
		\dot V(t,x,\delta x) \le -\alpha_1(V(t,x,\delta x))
	\end{equation}
	for a function $V(t,x,\delta x)\ge 0$ and  all relevant $(t,x,\delta x)$.  Integrating both sides  delivers for $t>0$ and with $c(t):=V(t,x(t),\delta x(t))$  that
	$$
	-\infty<-c(0)  \le c(t)-c(0)
	\le -\int_0^t \alpha_1(c(s)) ds
	\le -t\,  \alpha_1(c(t)).
	$$
	Since $c(t)\ge 0$  is monotonically decreasing by \eqref{genDEC}, $\lim_{t\to \infty} c(t)=c \ge 0$ and $c>0$ contradicts $-c(0) \le -t \alpha_1(c(t)) \le -t \alpha_1(c)$ for all  $t\ge 0$, we can conclude that $V(t,x(t),\delta x(t)) \to 0$ as $t\to \infty$.

	In the classical Lyapunov theory, one usually demands additionally to \eqref{genDEC}  that
	\begin{equation}
		\label{addlyacond}
		\alpha_2(\|x\|) \le V(t,x,\delta x) \le \alpha_3(\|x\|),
	\end{equation}
	and then $V(t,x(t),\delta x(t)) \to 0$ implies  $x(t)\to 0$ as $t \to \infty$; in general, one considers functions $V$ that do not depend on $\delta x$.  In contraction analysis, one usually demands additionally to \eqref{genDEC}  that
	\begin{equation}
		\label{addcontrond}
		\alpha_4(\|\delta x\|) \le V(t,x,\delta x) \le \alpha_5(\|\delta x\|),
	\end{equation}
	and then $V(t,x(t),\delta x(t)) \to 0$ implies  $\delta x(t)\to 0$ as $t \to \infty$.  Different choices of the functions $\alpha_i\in \cK_\infty$ and other variations in the conditions deliver different types of  stability.
	
	For simplicity of the exposition, we have considered flows on the manifold $\cM=\R^n$. In the classical Lyapunov theory one can also discuss flows on more general metric spaces $\cM$.
	In contraction analysis, one can consider flows of $\dot{x}=f(t,x)$ on  more general manifolds $\cM$ than $\R^n$, using $f\colon \R \times\cM \to T\cM$ with $T\cM=\{(x,v_x)\mid x\in \cM, v_x \in T_x\cM\}$ as the tangent bundle. One just needs some adequate notion of distance on the tangent spaces $T_x\cM$. Most commonly, one assumes that $\cM$ is a Riemannian manifold with (time-independent) metric $\langle v,w\rangle_x=v^TM(x)w$ (in local coordinates), where $v,w\in T_x \cM$ and $M(x)\in \cS^+_n$; note that the tangent space is isomorphic to $\R^n$. In practice, one often considers $M\colon \cM \to \cS^+_n$  as a variable of the problem and tries to determine it such that $V(x,\delta x)=\langle \delta x,\delta x\rangle_x$ fulfills \eqref{genDEC} and \eqref{addcontrond}; e.g.%
	\begin{equation*}
		\frac{\partial f^T}{\partial x}(t,x) M(x) + M(x) \frac{\partial f}{\partial x}(t,x)+
		\dot{M}(x)  \le -\beta M(x)\ \ \ \text{for}\ \ \ \beta \in\mathbb R^+.
	\end{equation*}
	Even more generally, one can consider Finsler metrics on the manifold $\cM$ as studied by Forni \& Sepulchre 2014 \cite{contr2014fornisepulchre}. A Finsler metric is given by a Finsler function on the tangent bundle $T\cM$, which equips each $T_x\cM$ with a (possibly asymmetric) Minkowskii norm instead of a scalar product as in the case of a Riemannian metric.  In this case, one must use the general formula \eqref{genLya2} in the decrease condition \eqref{genDEC}.
	%A \emph{Finsler-Lyapunov function}
	%
	%
	%The concept of a contraction metric can be generalized to a general Riemannian manifold $\cM$ as the phase space.
	%Then the Riemannian metric is defined by a scalar product on each tangent space $T_x \cM$; in the case of $\cM=\mathbb R^n$ above, it was given by $\langle v,w\rangle_x=v^TM(x)w$, where $v,w\in T_x \cM$; note that the tangent space is equivalent to $\mathbb R^n$.
	%
	%We will now generalize a Riemannian contraction metric to a Finsler-Lyapunov function. In short, a Riemannian contraction metric $M(x)$ defines a quadratic Finsler-Lyapunov function of the form $V(x,\delta x)=\delta x^TM(x)\delta x$.
	%
	%A Finsler metric is given by a Finsler function on the tangent bundle $T\cM$, which defines a Minkowskii norm on each tangent space. A {\bf Finsler-Lyapunov function}
	%$V\colon T  \cM \to \mathbb R_0^+$,
	%defined for every $(x,\delta x)\in T  \cM$, is such a Finsler function
	%% i.e. $$c_1F(x,\delta x)^p\le V(x,\delta x)\le c_2F(x,\delta x)^p$$ holds with a Finsler structure $F\colon T \cM\to \mathbb R_0^+$ and $p\ge 1$
	%with a contraction property of the following form
	%\begin{eqnarray}\frac{\partial V(x,\delta x)}{\partial x}f(x)+\frac{\partial V(x,\delta x)}{\partial \delta x}\frac{\partial f(x)}{\partial x}\delta x\le \alpha(V(x,\delta x))\qquad \mbox{ for all $\delta x\in T_x \cM$.}\label{cont2}
	%\end{eqnarray}
	%Here, different types of function $\alpha$ imply different types of incremental stability.

	Note that although contraction analysis is a more powerful tool than Lyapunov functions for the analysis of one attractor, Lyapunov functions are also very useful. In particular, the theory of Lyapunov functions for autonomous systems can be generalized to \emph{complete Lyapunov functions} that are defined on the whole state space, see Auslander 1964 \cite{complete1964auslander} and Conley 1978 \cite{complete1978conley}.  They are decreasing along all solution trajectories where possible, i.e.~on the complement of the chain-recurrent set, and their existence
	%fully characterize the long term behaviour
	has been established for general dynamical systems on separable metric spaces Hurley 1998 \cite{complete1998hurley}; see also  Hurley 1991--95 \cite{Hurley1991,Hurley1992,hurley95}, Akin 2010 \cite{akin2010topdyn}, and Patr\~{a}o 2011 \cite{complete2011patrao}.
	For flows on domains in  $\R^n$ (and easily extendable to manifolds), even the existence of a $C^\infty$ complete Lyapunov function has been established, see  Bernhard \& Suhr 2018 \cite{BeSu2018smmothConeLya} and Hafstein \& Suhr 2021 \cite{SuHa2020complya}. %, based on an even more general theorem for cone-fields.

	\subsection{Contraction analysis, incremental stability, and convergent systems}
	\label{sec:difftypesstab}
	There are three different notions which seek to compare the evolution of two trajectories: \emph{contraction analysis}, \emph{incremental stability},  and \emph{convergent systems}.
	
	Summarizing, contraction analysis seeks to formulate local, differential conditions for the contraction of two (adjacent) trajectories in a suitable metric, incremental stability requires the distance between any two solutions to decrease to zero as time goes to infinity, and convergent dynamics study the decrease of the distance to one specific solution $\bar{x}(t)$ as time goes to infinity.
	
	Thus, contraction analysis deals with a local criterion, incremental stability is interested in the evolution (contraction) of the distance between two solutions and convergent dynamics in convergence to one solution, so the limiting behavior. Often, in particular in autonomous systems and under additional conditions, contraction analysis is used to find conditions that imply incremental stability and then to show convergent dynamics, as exemplified for linear systems  in Section \ref{sec:linear}. However, in general non-autonomous systems, these three concepts of stability are independent and additional assumptions are needed for one to imply another. % neither notion implies the other ones.
	R\"uffer, van de Wouw \& Mueller 2013
	\cite{cont2013RWM} study the relation of convergent systems and incremental stability, while Fromion \& Scorletti 2005 \cite{cont2005FS} compare incremental stability with contraction analysis.
	We will later see different conditions, depending on the type of incremental stability and limiting behavior (e.g. equilibrium or periodic solution).
	
	Tran, R\"uffer \& Kellett 2019  \cite{disc2019TRK} consider the three notions \emph{incremental stability}, \emph{convergent dynamics} and \emph{contraction analysis} for discrete and time-varying systems. They compare the relations between the three properties and characterize them using Lyapunov functions.
	
	Let us compare convergent systems and incremental stability (see \cite{cont2013RWM}). %Let us denote the solution of the non-autonomous ODE $\dot{x}=f(t,x)$ with initial condition $x(t_0)=x_0$ by $\phi(t,t_0,x_0)$.
	A non-autonomous system is \emph{uniformly convergent} in a positively invariant set ${\mathcal X}\subset \mathbb{R}^n$ (see  Pliss 1966 \cite{cont1966Pliss}, who introduced (globally) convergent systems for time-periodic systems) if
	\begin{enumerate}
		\item all solutions $\phi(t,t_0,x_0)$ exist for all $t\in [t_0,\infty)$ for all initial conditions $(t_0,x_0)\in \mathbb R\times \mathcal X$,
		\item there exists a unique solution $\bar{x}(t)$ in $\mathcal X$, which is defined and bounded for $t\in\mathbb R$,
		\item the solution $\bar{x}(t)$ is uniformly asymptotically stable, i.e.~there exists a function $\beta \in {\mathcal KL}$ such that for all $(t_0,x_0)\in \mathbb R\times \mathcal X$ and $t\ge t_0$ we have
		$$		\|\phi(t,t_0,x_0)-\bar{x}(t)\|\le \beta(\|x_0-\bar{x}(t_0)\|,t-t_0).$$
	\end{enumerate}
	
	The system is \emph{incrementally asymptotically stable} (IAS) in a positively invariant set ${\mathcal X}\subset \mathbb R^n$ if  there exists a function $\beta \in {\mathcal KL}$ such that for all $\xi_1,\xi_2\in  \mathcal X$ and $t\ge t_0$ we have
	$$		\|\phi(t,t_0,\xi_1)-\phi(t,t_0,\xi_2)\|\le \beta(\|\xi_1-\xi_2\|,t-t_0).$$
	
	The system is \emph{globally incrementally stable} if   ${\mathcal X}= \mathbb R^n$. The system could have additional inputs (controls or disturbances from a closed set), see Angeli 2002  \cite{contr2002Angeli}.
	
	It is shown in \cite{cont2013RWM} that both notions are different in general: neither implies the other one. An example for a system that is uniformly convergent, but not globally IAS is a 2-dimensional system with solutions spiraling to a bounded solution, but at different angular velocities. An example for a system that is globally IAS, but not uniformly convergent is given by $\dot{x}=t-x$, since there does not exist a bounded solution. Note, however, that a reparameterization by subtracting a the limiting, unbounded solution $\tilde{x}(t)=t-1$ leads to a uniformly convergent system, i.e. $z=x-\tilde{x}(t)$ satisfies $\dot{z}=-z$.
	On the other hand, it is shown that both notions are equivalent if $\mathcal X$ is compact.
	
	Lastly, a characterization of globally uniformly convergent systems using Lyapunov functions is given in \cite[Th.~7]{cont2013RWM}, where the Lyapunov function $V\in C^1(\mathbb  R\times \mathbb R^n,\mathbb R_0^+)$ satisfies
	$$\alpha_1(\|x-\bar{x}(t)\|)\le V(t,x)\le \alpha_2(\|x-\bar{x}(t)\|)$$ and
	$$\dot{V}(t,x)\le -\alpha_3(\|x-\bar{x}(t)\|)$$
	with $\alpha_1,\alpha_2,\alpha_3\in {\mathcal K}_\infty$
	as well as $$V(t,0)\le c$$ for all $t\in\mathbb R$ and some constant $c>0$, which is equivalent to the boundedness of the solution $\bar{x}(t)$.
	
	\subsection{Historical overview of contraction analysis (equilibrium and time-periodic systems)}
	
	In this section, we discuss the earliest contributions to contraction analysis and partially follow the original notations.
	One of the earliest references is Trefftz 1926 \cite{Trefftz} who considers a second-order time-periodic ODE and defines stable solutions by the property that adjacent solutions (and their derivatives) converge to the stable solution (and its derivative) as time goes to infinity. It is shown that every stable solution converges to a periodic solution.
	Lewis 1949 \cite{contr1949lewis} uses geodesic distance in a Finsler space.
	Lewis 1951 \cite{cont1951Lewis}
	considers a time-independent metric and shows that if the Finsler metric $M_f(x,\dot{x})$ for a system $\dot{x}=f(t,x)$  satisfies
	$$\frac{\partial M_f}{\partial x}f+\frac{\partial M_f}{\partial \dot{x}}\frac{\partial f }{\partial x}\dot{x}\le -\beta<0$$
	in a region $R$ of the state space, then any two solutions must approach each other asymptotically.
	For a time-periodic system $\dot{x}=f(t,x)$, this implies that there exists a unique periodic solution to which all trajectories converge asymptotically, see also \cite{cont2005Jouffroy}.
	
	Reissig 1955 \cite{Reissig1955B} considers a second-order time-periodic ODE and defines the notion of {\it extreme stability} for the system, if the difference between any two solutions converges to zero as time goes to infinity. He shows that extreme stability is equivalent to the existence of a periodic solution to which all solutions converge.
	
	LaSalle 1957 \cite{LaSalle} defines an equivalence relation: two solutions are equivalent if the difference between them vanishes as time goes to infinity, leading to a partition of the phase space into equivalence classes of {\it contracting sets}.

	Seifert 1958 \cite{cont1958Seifert} considered a time-periodic ODE $\dot{x}=f(t,x)$ with $x\in \mathbb R^2$ and a compact, simply connected region $R\subset \mathbb R^2$. $R$ is assumed to have a certain type of positive invariance property. Seifert derives conditions on the boundary of $R$, such that the evolution of its length with respect to an appropriately chosen metric converges to $0$ and thus proves the existence, uniqueness and asymptotic stability of a periodic solution.
	
	Opial 1960 \cite{cont1960Opial} considers a two-dimensional, non-autonomous ODE and a  simply connected compact set $K\subset\mathbb R^2$, such that $\phi(t,t_0,x)\in K$ for all $t_0\in\mathbb R$, $t\ge t_0$, and $x\in K$. Then asymptotic stability of one solution is equivalent to incremental asymptotic stability of all solutions in the set, which he calls \emph{asymptotic stability of the system}.
	He measures the distance between two trajectories by the length of an appropriate curve between them.
	% $x_1(t)$ and $x_2(t)$,
	% namely
	%$$L(t)=\int_0^1 \left\|\frac{\partial x}{\partial s}(t,s)\right\|\,ds.$$
	%Setting $v(t,s)=\frac{\partial x}{\partial s}(t,s)$, he
	%studies the time evolution of $v$, namely
	%\begin{eqnarray}
	%\frac{dv}{dt}&=&\frac{\partial f}{\partial x}v
	%\end{eqnarray}
	If the time derivative of the quadratic form $v^TM(x)v$, where $M(x)\in {\mathcal S}_n^+$, satisfies an inequality for all $x\in K$, which in our notation is
	\begin{eqnarray}
		v^T\left(\frac{\partial f^T}{\partial x}M+M\frac{\partial f}{\partial x}+\dot{M}\right)v\le -\beta \|v\|^2
	\end{eqnarray}
	with $\beta\in\mathbb R^+$, then the system  is asymptotically stable, see \cite[Th.~2]{cont1960Opial}.
	Opial uses these results to compare the behavior of trajectories of different systems, using auxiliary systems (which is later called synchronization).

	Hartman 1961 \cite{cont1961Hartman}  considers a time-independent metric $M(x)=\Theta^T(x)\Theta(x)$  for the non-autonomous ODE $\dot{x}=f(t,x)$. He assumes%
	\begin{eqnarray*}
		\Theta^T\frac{\partial \Theta}{\partial x}f+M\frac{\partial f}{\partial x} < 0;
	\end{eqnarray*}
	note that in his notation $A<0$ for a possibly non-symmetric matrix means  $\frac{1}{2}(A^T+A)<0$ in our notation, i.e. the symmetric part
	\begin{eqnarray*}
		\left(	\Theta^T\frac{\partial \Theta}{\partial x}f+M\frac{\partial f}{\partial x}\right)^T+	\Theta^T\frac{\partial \Theta}{\partial x}f+M\frac{\partial f}{\partial x}&=&
		\frac{\partial f^T}{\partial x}M+M\frac{\partial f}{\partial x}+\dot{M},
	\end{eqnarray*}
	is negative definite. %, see \eqref{cont-M}.
	The conclusion is that the distance with respect to $M$ between any pair of solutions is decreasing.
	This condition asserts that if the system is autonomous and has an equilibrium (which can be shown by a stronger contraction condition), then it is globally asymptotically stable, see \cite[Lem.~1]{cont1961Hartman}.

	Hartman 1964 \cite{ODE1964hartman} introduces a criterion for global asymptotic stability of an equilibrium, which is assumed to exist, based on%
	\begin{eqnarray}
		f^T(x) \frac{\partial f}{\partial x}(x)f(x) \le 0\label{fc}
	\end{eqnarray}
	and another one based on
	\begin{eqnarray}
		v^T \frac{\partial f}{\partial x}(x)v \le 0\text{ for all }v\in\mathbb R^n\text{ with }v^Tf(x)=0.\label{sc}
	\end{eqnarray}
	He generalizes the Euclidean metric in \eqref{fc} to the case of a constant metric \cite[Ch.~XIV.10 and 11]{ODE1964hartman} and then further to a point-dependent Riemannian metric and shows implications of the existence of such a contraction metric in \cite[Ch.~XIV.12 and 13]{ODE1964hartman}. Further, the criterion \eqref{sc} is extended to general contraction metrics in \cite[Ch.~XIV.14 to 16]{ODE1964hartman}, showing contraction between trajectories and deducing a globally stable equilibrium -- note that the criterion is also related to the existence of limit cycles if there are no equilibria,  % which is also mentioned in the book,
	see Section \ref{sec:period}.

	Krasovski\u{i} considers a non-autonomous ODE with solution $x(t)\equiv 0$ and a constant metric $M$. He derives a criterion that shows global asymptotic stability of an equilibrium, see Krasovski\u{i} 1963 \cite[Th. 21.1]{stability1963krasovskii}, based on Krasovski\u{i} 1954 \cite{krasovskii1954} and Krasovski\u{i} 1957 \cite{krasovskii1957}. This theorem is also discussed in Hahn 1967 \cite[Th.~55.5]{stability1967hahn} in the autonomous case; the proof uses the Lyapunov function $v(x)=x^TMx$.

	Demidovi\v{c} 1961 \cite{cont1961Demidovic} considers a non-autonomous ODE and uses a constant (time- and space-independent) contraction metric to show that all solutions are exponentially stable. He also shows that the system then is dissipative, i.e. every solution starting at $x$ is in a fixed region of the state space after some time $T(x)$.

	Demidovi\v{c} 1967 \cite{cont1967Demidovic} defined and studied convergent systems for general non-autonomous ODEs.
	He showed that if   \begin{eqnarray}J(t,x) = \frac{1}{2}\left( \frac{\partial f^T}{\partial x}(t,x) M+M\frac{\partial f}{\partial x}(t,x) \right)\label{J-demi}
	\end{eqnarray} is negative definite uniformly in $(t,x)\in\mathbb R\times \mathbb R^n$, where $M\in {\mathcal S}^+_n$, then the distance between any two solutions decreases exponentially. If $|f(t,x)|\le c<\infty$ then, in particular, no finite escape times exist and all solutions are globally uniformly exponentially stable. The proof considers the function $\frac{1}{2}(x_1(t)-x_2(t))^TM(x_1(t)-x_2(t))$ along two solutions, and estimates the derivative using the mean value theorem.
	If $f(t,0)=0$, then this is Krasovski\u{i}'s
	stability theorem, see Krasovski\u{i}'s 1963 \cite{stability1963krasovskii}.
	This result can be applied in case that the Jacobian of $f$ is in the convex hull of matrices $A_1,\ldots,A_k$ and there is a solution of the linear matrix inequalities $A_i^TM+MA_i<0$.

	Datko 1966 \cite{datko1966} assumes the existence of an asymptotically stable equilibrium and a constant contraction metric $M\in {\mathcal S}^+_n$ such that
	$$\frac{\partial f^T}{\partial x}(x)M+M\frac{\partial f}{\partial x}(x)\le 0.$$
	
	\vspace{0.3cm}

	Lohmiller \& Slotine 1998 \cite{contr1998lohmillerslotine} is one of the key references from the 1990s in contraction analysis, which has been cited by many authors. They consider the non-autonomous differential equation
	\begin{eqnarray}
		\dot{x} = f(t,x)\label{ODE}
	\end{eqnarray}
	where $x\in \mathbb R^n$.
	
	The main idea is to study the time evolution of the distance between two adjacent solutions of \eqref{ODE}. This can either be done by two solutions, starting a small distance apart, or by infinitesimal displacements, using the virtual displacement $\delta x$. The dynamics for these virtual displacements are governed by%
	\begin{eqnarray}
		\dot{\delta x} = \frac{\partial f}{\partial x}(t,x) \delta x. \label{variation}
	\end{eqnarray}
	It is shown that if $\frac{\partial f}{\partial x}(t,x)$ is uniformly negative definite, i.e.,
	$$\exists b>0\quad \forall x\in \mathbb R^n \quad \forall t\ge 0
	\colon\quad \frac{1}{2}\left(\frac{\partial f}{\partial x}(t,x)
	+\frac{\partial f^T}{\partial x}(t,x)\right)\le -b I,$$
	then any (second) trajectory which starts in a ball of constant radius centered about a given (first) trajectory, and contained in the region where the above condition holds (the contraction region), remains in that ball and converges exponentially to the first trajectory.
	This result is a strengthened version of Krasovski\u{i}'s  theorem on global asymptotic convergence, see Krasovski\u{i} 1963 \cite{stability1963krasovskii}.
	While these results rely on the use of the Euclidean metric, another ingredient is a coordinate change or a different Riemannian metric, see Section \ref{sec:diffformcontr} or \eqref{zformF}.

	A region in the state space is called \emph{contraction region} if, see \eqref{refF} for the formula for $F$,
	\begin{eqnarray}
		F(t,x)&\le&-b I\text{ or equivalently}\label{contr00}\\
		\frac{\partial f^T}{\partial x}(t,x)M(t,x)+M(t,x)\frac{\partial f}{\partial x}(t,x) +\dot{M}(t,x)
		&\le&-2b M(t,x).\label{contraction}
	\end{eqnarray}
	Note that the notation \eqref{contr00} means that the symmetric part of $F(t,x)+bI$ is negative semidefinite; for the equivalence, see \eqref{contraction0}.
	
	Any first trajectory which remains in the contraction region for all positive times, and any second trajectory which starts in a ball of constant radius with respect to the metric $M(t,x)$ remains in that ball and converges exponentially to the first trajectory.
	Semi-contraction is defined by $F(t,x)\le 0$ or equivalently the left-hand side of \eqref{contraction} being negative semi-definite.
	
	Some consequences include that (i) a convex contraction region contains at most one equilibrium, (ii) in an autonomous, globally contracting system all trajectories converge exponentially to a unique equilibrium, (iii) a time-periodic, globally contracting system tends exponentially to a periodic solution.

	The paper also introduces semi-contraction differently by assuming that $F$ is uniformly negative definite, but $M$ is only positive semi-definite with some principal directions corresponding to uniformly positive eigenvalues of $M$. This results in exponential convergence to zero of the components of $\delta x$ on the subspace spanned by the principal directions mentioned above.
	
	Aminzare \& Sontag 2014 \cite{cont2014AS} consider a metric that is independent of both $t$ and $x$, defined by a positive definite (constant) matrix, and also generalize the conditions to infinite-dimensional systems given by reaction-diffusion systems.
	
	Margaliot, Sontag \& Tuller 2016
	\cite{cont2016MST} relax the  classical contraction condition by considering contraction after a small transition time.

	Botner, Zarai, Margaliot \& Gr\"une 2017
	\cite{contr2017BZMG} compare the solutions of two systems $\dot{x}=f(t,x)$ and $\dot{y}=g(t,y)$ and assume that the first one is contracting. They derive bounds on the difference between solutions $d(t)=x(t)-y(t)$. This is applied to the case when the $y$-system has ``simple'' solutions, e.g. if they are known explicitly, or when $\|f-g\|$ admits a simple bound as this term appears in the estimate.

	\subsection{Reviews and tutorials}
	
	The review Jouffroy 2005 \cite{cont2005Jouffroy} puts the definition from Lohmiller \& Slotine 1998 \cite{contr1998lohmillerslotine} into historical context and cites earlier discussions of contraction criteria, including some of the ones mentioned above.
	The review ``Convergent Dynamics'' Pavlov et al. 2004 \cite{rev2004PPWN} introduces the contributions of Demidovi\v{c} and puts them into perspective with other results, the authors also include the main ideas and some proofs.
	Sontag 2010 \cite{Son2010} contains an introduction to contraction analysis and proves some basic results.
	%Advantage of stability of {\it all} solutions versus stability of a specific solution that (a) specific solution does not need to be found and (b) synchronization
	Simpson-Porco \& Bullo 2014 \cite{man2014SB} present contraction analysis on Riemannian manifolds, rigorously using intrinsic geometric concepts.
	
	There are several books about contraction metrics such as  Lohmiller \& Slotine 2000 \cite{cont2000LS}, an overview and open questions Aminzare \& Sontag 2014 \cite{cont2014AS}, and
	a tutorial on incremental stability and contraction Jouffroy \& Fossen 2010 \cite{cont2010JF}. The latter includes a discussion of the conditions
	\begin{enumerate}
		\item[(i)] $\frac{\partial f^T}{\partial x}M+M\frac{\partial f}{\partial x}+\dot{M}\le -\beta I$
		\item[(ii)] $\frac{\partial f^T}{\partial x}M+M\frac{\partial f}{\partial x}+\dot{M}\le -\beta M$
	\end{enumerate}
	As $M(t,x)$ may become unbounded for $t\to\infty$, the condition (ii)  guarantees exponential convergence, while the condition (i) does not.
	
	%Partial contraction compares two systems of the form
	%\begin{eqnarray}
	%\dot{x}&=&f(x,x,t)\\
	%\dot{y}&=&f(y,x,t)
	%\end{eqnarray}
	%If the second system is contracting with respect to $y$ and a particular solution of the second system verifies a smooth specific property, then all trajectories of the original (first) system satisfy this property exponentially and the first system is called partially contracting., see also \cite{cont2003Slotine,cont2005WS}.
	
	\emph{Partial contraction analysis}, see Slotine \& Wang 2005 \cite{Slotine2005} or the review Jouffroy \& Fossen 2010 \cite{cont2010JF}, constructs an auxiliary (virtual) system, which is contracting with respect to an additional, auxiliary variable. If a particular solution of this auxiliary system satisfies a specific property, then all trajectories of the original system satisfy this property exponentially, and the original system is called partially contracting, see also \cite{cont2003Slotine,cont2005WS}.

	\subsection{Converse theorems for equilibria}
	\label{sec:converse}
	
	The existence of a contraction metric implies (under certain conditions) that solutions converge exponentially to a unique equilibrium. Here we deal with the converse, i.e.~whether, given an exponentially stable equilibrium, a contraction metric exists. We have already discussed this in the case of autonomous linear systems in Section \ref{sec:aut-lin}.

	Lohmiller \&  Slotine 1998 \cite{contr1998lohmillerslotine} assume that a non-autonomous system is exponentially convergent with rate $\beta/2$. They then define a metric $M(t,x(t))$ for each trajectory as the solution of
	\begin{eqnarray}\dot{M}&=&-\beta M-M\frac{\partial f}{\partial x}-\frac{\partial f^T}{\partial x}M\\
		M(0,x(0))&=&kI
	\end{eqnarray}and show that $M$ is uniformly positive definite with $M\ge I$. However, this metric $M$ depends on time and may become unbounded as $t\to\infty$.
	
	Giesl 2005 \cite{contr2015giesl} considers an autonomous equation $\dot{x}=f(x)$ with an exponentially stable equilibrium and proves  the existence of three
	contraction metrics: the first is defined on a compact subset of the basin of attraction, and satisfies a certain inequality (Theorem 4.1); the second is defined  on the whole basin of attraction, recovers the exponential attraction rate arbitrarily well, but is only continuous and orbitally continuously differentiable  (Theorem 4.2); the third is defined  on the whole basin of attraction, is as smooth as $f$ and satisfies a matrix-valued PDE (Theorem 4.4), which is useful for constructive methods, see Section \ref{sec:collo}.
	%
	%Note that we consider contraction in directions $\bv$ perpendicular to $f(\bx)$ with respect to the  metric $M$, i.e. $\bv^TM(\bx)f(\bx)=0$. One could alternatively consider directions perpendicular to $f(\bx)$ with respect to the Euclidean metric, i.e. $\bv^Tf(\bx)=0$, but then the function $L_M$ needs to reflect this, see \cite{contGiesl2021,cont1988BL}.
	
	%
	%Converse theorem for periodic orbits
	%\cite{contr2019bGiesl}
	%
	%Characterized by PDE: \cite{contGiesl2021}
	
	\section{Extensions of contraction analysis}
	\label{sec:ext}
	
	In this section, we present different extensions of contraction analysis. We start with contraction analysis for periodic orbits, which only requires contraction in directions transversal to the flow.  Then we consider generalizations of contraction analysis to different classes of systems.
	Note however, that there are numerous applications of contraction analysis to other types of systems studied in the literature and we only highlight a selection  in this review.
	We will also discuss some additional types, for which  numerical methods have been developed, in Section \ref{sec:num}; in particular, we will present \emph{switched systems}, \emph{differential-algebraic systems}, and \emph{reaction-diffusion equations} together with related \emph{compartmental ODEs}.
	
	\subsection{Contraction analysis for periodic orbits}
	\label{sec:period}

	If a system has an asymptotically stable periodic orbit, rather than an asymptotically stable equilibrium, then contraction in direction of the flow $f(x)$ cannot occur. Hence, the contraction condition
	$$\frac{\partial f^T}{\partial x}(x)M(x)+M(x)\frac{\partial f}{\partial x}(x)+\dot{M}(x)
	\le -\beta M(x),$$
	which is equivalent to
	\begin{eqnarray}
		v^T\left(\frac{\partial f^T}{\partial x}(x)M(x)+M(x)\frac{\partial f}{\partial x}(x)+\dot{M}(x)\right)v
		\le -\beta v^TM(x)v\label{cont-per1}
	\end{eqnarray} for all $v\in\mathbb R^n$
	is relaxed to \eqref{cont-per1} for all $v\in\mathbb R^n$ with $v^TM(x)f(x)=0$. This restricts the contraction condition of adjacent solutions to difference vectors $v$ that are perpendicular to the flow at $x$ with respect to the metric $M$. One can, similarly to Section
	\ref{sec:linear}, show that the distance between adjacent trajectories $\phi(t,x)$ and $\phi(t,y)$ decreases exponentially, however, $y-x$ needs to be perpendicular to $f(x)$  in the metric $M$ and the time needs to be synchronized such that the difference vector is perpendicular at each time $t$. The stability with time synchronization of trajectories is called \emph{Zhukovski stability} and is equivalent to \emph{orbital stability}, see Leonov, Ponomarenko \& Smirnova 1995 \cite{LPS1995} and
	Leonov 2006 \cite{per2006Leonov}.%
	
	One can also synchronize the times such that the difference vector $v$ is perpendicular to $f(x)$ with respect to the Euclidean metric -- this has advantages when constructing the metric $M$. Then the condition becomes
	\begin{eqnarray}\lefteqn{
		v^T\left(V^T(x)M(x)+M(x)V(x)+\dot{M}(x)\right)v}\nonumber\\
		&\le&-\beta v^TM(x)v\text{\quad for all }v^Tf(x)=0\label{peri0}\\
		\text{ where }
		V(x)&:=&\frac{\partial f}{\partial x}(x)-\frac{f(x)f^T(x)(\frac{\partial f^T}{\partial x}(x)+\frac{\partial f}{\partial x}(x))}{\|f(x)\|^2}. \label{peri1c}
	\end{eqnarray}
	
		\begin{figure}[ht]
		\centering
		\includegraphics[width=0.9\textwidth]{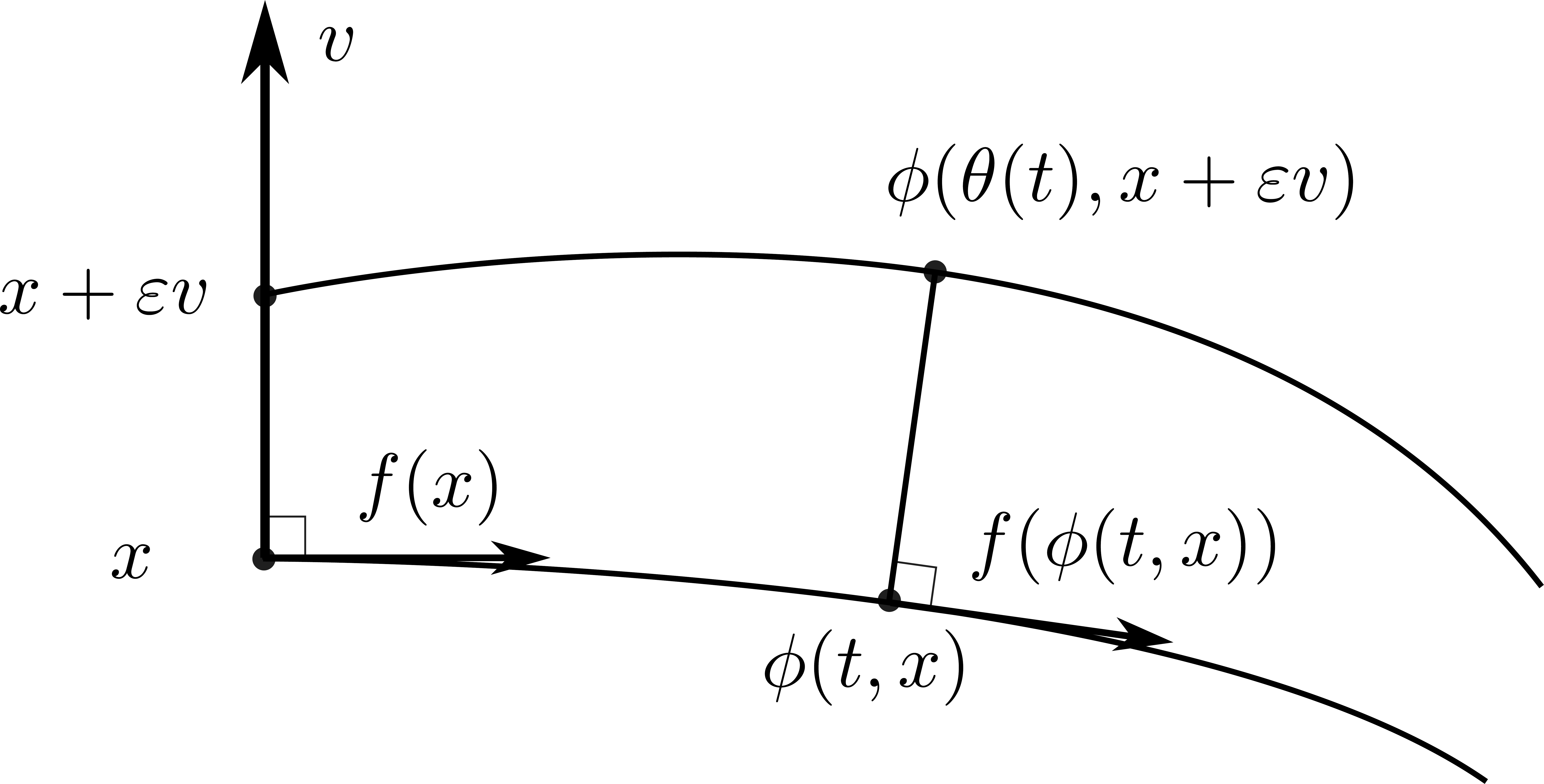}
		\caption{Solutions starting at $x$ and $x+\varepsilon v$, where $v$ is perpendicular to the flow $f(x)$ at $x$.  The time of the solution starting at $x+\varepsilon v$ is then synchronized through $\theta(t)$ such that the difference vector $\phi(\theta(t),x+\varepsilon v)-\phi(t,x)$ is
			perpendicular to the flow $f(\phi(t,x))$ at $\phi(t,x)$ for all times $t$. 
		}
		\label{fig4}
	\end{figure}

	Let us give a detailed explanation of this condition: consider two adjacent solutions at $\phi(t,x)$ and $\phi(t,x+\varepsilon v)$ for small $\varepsilon>0$, where $v^T f(x)=0$, so $v$ is perpendicular to $f(x)$.
	% We will give a heuristic argument t that $$v^T\left(\dot{M}(x)+V(x)^TM(x)+M(x)V(x)\right)v<0$$ for all $v^Tf(x)$, see the left-hand side of \eqref{peri1}, if the distance between solutions through $\phi(t,x)$ and $\phi(t,x+\varepsilon v)$ with respect to the metric $M(x)$ decreases.
	To measure the distance, we synchronize the time, i.e.~we consider the two solution trajectories $\phi(t,x)$ and $\phi(\theta(t),x+\varepsilon v)$, such that the difference vector  $\phi(\theta(t),x+\varepsilon v)-\phi(t,x)$ is perpendicular to $f(\phi(t,x))$ with respect to the Euclidean metric; for a schematic depiction see Figure \ref{fig4}. In particular, we define $\theta(t)$ such that $\theta(0)=0$ and
	\begin{eqnarray}\left(\phi(\theta(t),x+\varepsilon v)-\phi(t,x)\right)^T f(\phi(t,x))=0\text{ for all }t\ge 0.\label{per0}
	\end{eqnarray}
This is possible due to the implicit function theorem if $\varepsilon>0$ is sufficiently small and, in the first instance, also $t\ge 0$ is sufficiently small; it will follow for all $t\ge 0$ from the contraction condition.

	The implicit function theorem shows that
	\begin{eqnarray}
		\dot{\theta}(0) = \frac{\|f(x)\|^2-\varepsilon v^T \frac{\partial f}{\partial x}(x) f(x)}{f(x+\varepsilon v)^Tf(x)}\
		\approx \ 1-\varepsilon \frac{v^T (\frac{\partial f^T}{\partial x}(x)+\frac{\partial f}{\partial x}(x))f(x)}{\|f(x)\|^2+\varepsilon v^T \frac{\partial f^T}{\partial x}(x)f(x)}\label{thetadot}
	\end{eqnarray}
	for small $\varepsilon>0$.

	Now we consider the squared distance between the trajectories with respect to the Riemannian metric
	\begin{eqnarray}d(t) = \left(\phi(\theta(t),x+\varepsilon v)-\phi(t,x)\right)^TM(\phi(t,x))\left(\phi(\theta(t),x+\varepsilon v)-\phi(t,x)\right)\label{d0}
	\end{eqnarray} and take the derivative. We obtain for small $\varepsilon>0$, using Taylor expansion and keeping terms only up to order $\varepsilon^2$
	\begin{eqnarray*}
		\frac{d}{dt}d(t)\bigg|_{t=0}
		&=&\left(\dot{\theta}(0)f(x+\varepsilon v)-f(x)\right)^TM(x)\varepsilon v
		+\varepsilon^2 v^T\dot{M}(x)v \\	&&
		+\varepsilon v^T
		M(x)\left(\dot{\theta}(0)f(x+\varepsilon v)-f(x)\right)\\
		&\approx&\varepsilon(\dot{\theta}(0)-1)
		[f^T(x)M(x)v+v^TM(x)f(x)]\\
		&&
		+\varepsilon^2 \dot{\theta}(0)\left[\left(\frac{\partial f}{\partial x }(x)v\right)^TM(x)v
		+v^TM(x)\frac{\partial f}{\partial x }(x)v\right]
		+\varepsilon^2 v^T\dot{M}(x)v	\end{eqnarray*}
				\begin{eqnarray*}
		&\approx&\varepsilon^2\bigg[
		-\frac{v^T (\frac{\partial f^T}{\partial x }(x)+\frac{\partial f}{\partial x }(x))f(x)}{\|f(x)\|^2}[f(x)^TM(x)v+v^TM(x)f(x)]\\	&&
		+\left(\frac{\partial f}{\partial x }(x)v\right)^TM(x)v
		+v^TM(x)\frac{\partial f}{\partial x }(x)v
		+v^T\dot{M}(x)v
		\bigg]
		\text{ by }\eqref{thetadot}\\
		&=&\varepsilon^2v^T\left[V^T(x)M(x)+M(x)V(x)
		%-\frac{D\bff(\bx)\bff(\bx)\cdot \bff(\bx)^TM(\bx)\bv}{\|\bff(\bx)\|^2}
		%+D\bff(\bx)^TM(\bx)
		%+M(\bx)D\bff(\bx)
		+\dot{M}(x)\right]v\\
		&\le&-\beta  d(t),
	\end{eqnarray*}
	by \eqref{peri0}, which shows exponential convergence to zero.
	
	Finally, note that conditions such as
	\begin{eqnarray}
		v^T\left(\frac{\partial f^T}{\partial x}(x)M(x)+M(x)\frac{\partial f}{\partial x}(x)+\dot{M}(x)\right)v
		\le 0\text{\quad for all }v^Tf(x)=0\label{con1}
	\end{eqnarray}
	can also be expressed, under certain conditions, by a matrix inequality of the form
	\begin{eqnarray}
		\frac{\partial f^T}{\partial x}(x)M(x)+M(x)\frac{\partial f}{\partial x}(x)+\dot{M}(x)-cf(x)f^T(x)
		\le 0\label{con2}
	\end{eqnarray}
	for sufficiently large $c>0$. The idea is that \eqref{con2} implies \eqref{con1}, since $v^Tf(x)f^T(x)v=0$, while \eqref{con1} implies \eqref{con2} for sufficiently large $c>0$.
	
	\subsection{Historical overview of contraction analysis (periodic orbit)}
	
	Borg 1960 \cite{cont1960borg} considers an autonomous ODE and studies a contraction condition transversal to the flow, using the Euclidean metric, e.g.%
	\begin{eqnarray}
		v^T\left(\frac{\partial f^T}{\partial x}(x)+\frac{\partial f}{\partial x}(x)\right)v \le -\beta \|v\|^2\ \ \ \text{for all $v\in\mathbb R^n$ with  $v^T f(x)=0$.}\label{cond}
	\end{eqnarray}
	He shows that if this condition is satisfied in a positively invariant set which contains no equilibrium, then all solutions starting in this set converge to a unique periodic orbit in this set.  He uses a synchronization of the time of two adjacent trajectories such that their difference is perpendicular to the flow through the implicit function theorem in order to use \eqref{cond}. In an earlier version, Borg 1953 \cite{borg1953} formulates the condition in terms of $\lim_{\rho\to 0}\frac{v^T (f(x+\rho v )-f(x))}{\rho}<0$ for all $v\in\mathbb R^n\setminus \{0\}$ with  $v^T f(x)=0$ with the implication that all solutions converge to an equilibrium or a periodic orbit. Borg's criterion was applied in Sherman 1963 \cite{Sherman1963} to a 3-dimensional autonomous ODE arising from a nuclear spin generator.%
	
	Stenstroem 1962 \cite{contr1962stenstroem} generalizes these results to equations on a Riemannian manifold. He drops the assumption that the set does not contain an equilibrium and shows that then either solutions converge to an equilibrium or to a periodic orbit.%
	
	Hartman \& Olech 1962 \cite{contr1962hartmanolech} consider this condition, using a point-dependent metric on a Riemannian manifold and, moreover,  derive a sufficient condition based on the sum of the two largest eigenvalues. This idea can be generalized to using metrics for estimating the dimension of attractors, see Section \ref{sec:dimension}. They also consider non-autonomous equations.
	
	Leonov 1987 \cite{Leonov1987} generalizes these results to point-dependent metrics $M(x)$ and replaces the synchronization condition $v^Tf(x)=0$ by the more general condition $v^Tq(x)=0$, where $q(x)=f(x)$ is a special case.
	
	The Markus-Yamabe conjecture is that an equilibrium at $0$ of an autonomous ODE is globally stable if all eigenvalues of $\frac{\partial f}{\partial x}(x)$ have a negative  real part for all $x\in\mathbb R^n$; this goes back to Aizerman 1949 \cite{aizerman1949}. The conjecture has later been shown to be true for dimensions $n\le 2$, and false for dimensions $n\ge 3$. Markus \& Yamabe 1960 \cite{MY1960} have derived a criterion based on the negativity of the eigenvalues of the symmetric Jacobian $\frac{\partial f^T}{\partial x}(x)+\frac{\partial f}{\partial x}(x)$, which is based on a more general result about contraction metrics on Riemannian manifolds. Hartman 1961 \cite{cont1961Hartman} generalizes this result to Riemannian manifolds and Furi, Martelli \& O'Neill 2009 \cite{Furi2009} further generalize the result by assuming less smoothness on $f$.%
	
	Cronin 1980 \cite {Cronin1980} presents a contraction criterion requiring that $\frac{ \partial f}{\partial x}(x)$ has $n-1$ eigenvalues with negative real parts, and the eigenvector corresponding to the remaining eigenvalue is not perpendicular to $f(x)$.

	The criterion \begin{eqnarray}
		V^T(x)M(x)+M(x)V(x)+\dot{M}(x)\le -2\nu I
	\end{eqnarray}
	with $V$ given by \eqref{peri1c}  was introduced by Biochenko \& Leonov 1988 \cite{cont1988BL}.
	
	In Leonov 1990
	\cite{Leonov1990}, a point-dependent contraction metric with contraction in directions $v$ with  $v^Tq(x)=0$, where $q^T(x)f(x)\not =0$, is considered; this generalizes the choice $q(x)=f(x)$. Then trajectories approach each other as $t\to \infty$ and all converge to an equilibrium, if one exists.

	Leonov 2006
	\cite{per2006Leonov} has defined a moving Poincar\'e section to synchronize the time of trajectories. The paper links the definitions of Poincar\'e (or orbital stability) with Zhukovski stability, see also \cite{cont1998Leonov,leonov2001}, i.e. stability of solutions after reparameterization of time.  Here, the synchronization is such that the difference vector is perpendicular to $f$ with respect to the Euclidean metric. He considers Zhukovski stability of a general solution; in the case of a periodic solution, these results are the Andronov-Vitt theorem.
	The reparameterization or synchronization of the time of adjacent trajectories is used to show that the existence of a contraction metric implies the existence  a unique, exponentially stable periodic orbit to which all trajectories converge, see also Kravchuk,  Leonov \& Ponomarenko 1992 \cite{contr1992kravchukleonovponomarenko}, Yang 2001 \cite{Yang} or Manchester \& Slotine 2014 \cite{contr2014manchesterslotine}.

	\subsection{Converse theorems for periodic orbits}
	
	Contraction metrics, where contraction only occurs in an $(n-1)$-dimensional hyperplane, which is transversal to the flow $f$, imply (under certain conditions) that solutions exponentially converge to a unique periodic orbit. The converse question, namely the existence of such a contraction metric on the basin of attraction of a given exponentially stable periodic orbit, has been considered by several authors. A key ingredient, as for the sufficiency, is the synchronization of the time of adjacent trajectories.
	
	Converse theorems to prove the existence of a Riemannian contraction metric for a periodic orbit go back to Boichenko \& Leonov 1988 \cite[Th.~1]{cont1988BL}, constructing a metric $M(t)$ depending on time. A local version is presented in Hauser \& Chung 1994 \cite{per1994HC}, proving a converse theorem to construct a quadratic Lyapunov function for an exponentially stable periodic orbit. They use a local coordinate system around the periodic orbit, separating the tangential and transverse dynamics, and thus construct a local Lyapunov function -- valid in a (possibly small) neighborhood of the periodic orbit.
	
	A global converse theorem on a compact set is shown in Manchester \& Slotine 2014 \cite{contr2014manchesterslotine}, where also the robustness to parameters is discussed, using the construction in Leonov 2006 \cite{per2006Leonov}.
	
	In Giesl 2019 \cite{contr2019bGiesl},  the existence of a contraction metric in the entire basin of attraction is shown which recovers the rate of contraction arbitrarily well; here, the synchronization is such that the difference vector between two adjacent solutions is perpendicular to $f$ with respect to the Riemannian metric.
	
	Giesl 2021  \cite{contGiesl2021} uses the synchronization with respect to the Euclidean metric as in \cite{contr2014manchesterslotine}, and  shows the existence of a contraction metric in the entire basin of attraction $A(\Omega)$ of the periodic orbit $\Omega$ of the ODE $\dot{x}=f(x)$ with $f\in C^s(\mathbb R^n,\mathbb R^n)$. The contraction metric $M\in C^{s-1}( A(\Omega),{\mathcal S}_n^+)$ is the unique solution of the linear PDE \eqref{peri1} with any fixed right-hand side $B\in C^{s-1}( A(\Omega),{\mathcal S}_n^+)$, and \eqref{peri2} with any point $x_0\in A(\Omega)$ and $c_0\in \mathbb R^+$:
	\begin{eqnarray}V(x)^TM(x)+M(x)V(x)+\dot{M}(x)&=&-P_xB(x)P_x \text{ and}\label{peri1}\\
		f^T(x_0)M(x_0)f(x_0)&=&c_0\|f(x_0)\|^4 \text{, where}\label{peri2}\\
		P_x&=&I-\frac{f(x)f^T(x)}{\|f(x)\|^2}\label{peri3}
	\end{eqnarray}
	and  $V$ is defined as in \eqref{peri1c}.
	This can be used to approximate $M$ by solving the PDE numerically. Note that \eqref{peri1} ensures that the left-hand side is negative definite in directions perpendicular to $f$, as $P_x$, see \eqref{peri3}, projects onto the $(n-1)$-dimensional hyperplane perpendicular to $f$. Fixing $M$ in $f$-direction at one point $x_0$ in \eqref{peri2} is sufficient together with \eqref{peri1} to ensure that $M$ is positive definite in $A(\Omega)$.
	
	\vspace{0.3cm}
	
	In the following sections, we present some extensions and aspects of contraction metrics that have been studied in the literature. We apologize in advance for the contributions we missed or could not include because of space constraints.

	\subsection{Uncertain/stochastic systems}
	%Stochastic check \cite{pham2013}
	%
	%and \cite{BoSl2019contSto}
	%
	%\cite{stoch2019ISM} -- but this is Lyapunov theory?!  JA, UND SCHLECHT DAZU :)
	%
	%\cite{stoch2018Flynn} CHECK?  EIGENTLICH NICHT CONTRACTION IM UNSEREM SINNE
	%
	%\cite{stoch2009PTS} CHECK?
	%SIGGI, hatten wir nicht noch andere???
	
	In Pham, Tabareau \& Slotine 2009 \cite{stoch2009PTS}, the incremental stability of It\^{o} stochastic dynamical systems is studied; in particular \emph{stochastic differential equations} (SDEs)  of the form
	\begin{equation}
		\label{SDE}
		dX(t)=f(t,X(t))dt+G(t,X(t))dW(t),
	\end{equation}
	where $f\colon\R^+_0\times \R^n\to\R^n$, $G\colon\R^+_0\times\R^{n\times Q}\to \R^n$ and $W\colon\R^+_0\to\R^Q$ is a $Q$-dimensional Wiener process.
	It is assumed that $f(t,\cdot)$ and $G(t,\cdot)$ are globally Lipschitz uniformly in $t\ge0$ and then strong solutions
	$$
	X(t)=X^0+\int_{t_0}^tf(\tau,X(\tau))d\tau+\int_{t_0}^tG(\tau,X(\tau))dW(\tau),
	$$
	that are unique and continuous (a.s.) exist.
	Here $X(t_0)=X^0$ in $\R^n$ is the initial distribution and
	the underlying probability space and the filtration satisfy the \emph{usual conditions}, cf.~e.g.~\cite[\S21]{sdestab2012khaminskii} or \cite[\S2.3]{sdestab2008mao}. The second integral is interpreted in the It\^{o} sense; note that this also includes the Stratonovich interpretation of the SDE
		with modified drift coefficient.
	%for a slightly modified SDE.
	
	Various stability properties of the zero solution of \eqref{SDE} with $f(t,0)=0$ and  $G(t,0)=0$ for all $t\ge0$ have been  studied using Lyapunov functions in the literature, both theoretically \cite{sdestab2012khaminskii,sdestab2008mao,GuHa2018theo} and numerically \cite{HGGS2018localLya,BHGSG2019stobas,CaGr2003stoZu}. {Caraballo, Kloeden \& Schmalfu\ss \, 2004 \cite{Caraballo2004} proved the existence of exponentially stable non-trivial solutions to stochastic semilinear partial differential equations generated by random variables chosen as initial values. Using the pullback method, they establish a link between the local analysis of the temporal behavior of stochastic partial differential equations and global analysis in random dynamical systems, see e.g.~Arnold \cite{Arno2002rdsbook}.}

	In Pham, Tabareau \& Slotine 2009 \cite{stoch2009PTS}, the incremental stability of the system
	\begin{equation}
		\label{SDE2}
		d\begin{pmatrix}
			X_1(t) \\
			X_2(t) \\
		\end{pmatrix} = \begin{pmatrix}
			f_1(t,X_1(t)) \\
			f_2(t,X_2(t)) \\
		\end{pmatrix}dt+ \begin{pmatrix}
			G(t,X_1(t)) & 0 \\
			0 & G(t,X_2(t)) \\
		\end{pmatrix}dW(t)
	\end{equation}
	is studied, where $f_1,f_2$ and $G_1,G_2$ are of the same dimensions as $f$ and $G$ in \eqref{SDE}, respectively,  and $W$ is a $2Q$-dimensional Wiener process.
	Note that \eqref{SDE2} can be interpreted as two different SDEs of the form \eqref{SDE} with independent Wiener processes.

	First, the case $f_1=f_2=f$ and $G_1=G_2=G$ is studied, under the assumption that the deterministic system $\dot{x}=f(t,x)$ is uniformly contracting with respect to a time-varying metric $M\colon \R^+_0\to \cS^+_n$  and that $\tr\left(G^T(t,x)M(t)G(t,x)\right)$ is uniformly bounded.
	This means that there are constants $a,b,c>0$ such that
	$$
	M(t)\ge a I\ \ \ \text{and}\ \ \ \frac{\partial f^T}{\partial x}(t,x)M(t)+  M(t)\frac{\partial f}{\partial x}(t,x)+\dot{M}(t) \le -2b M(t)
	%M(t)\ge \alpha I\ \ \ \text{and}\ \ \ M(t)\frac{\partial f}{\partial x}(t,x)+\frac{\partial f}{\partial x}(t,x)^TM(t)+\dot{M}(t) \le -2\beta M(t)
	$$
	and
	$$
	\tr\left(G^T(t,x)M(t)G(t,x)\right)\le c < \infty
	$$
	for all $(t,x)\in \R^+_0\times \R^n$, and a system fulfilling these conditions is said to be \emph{stochastically contracting} in the metric $M$. % with rate $\beta$ and bound $C$.
	It is proved by using the Lyapunov function $V(t,x_1,x_2)=(x_1-x_2)^TM(t)(x_1-x_2)$, that a solution $(X_1(t),X_2(t))^T$
	with initial distribution %$(X_1(0),X_2(0))=
	$(X_1^0,X_2^0)$ at time $t=0$ fulfills
	\begin{equation}
		\label{stocontr1}
		\mathbb {E}(\|X_1(t)-X_2(t)\|^2) \le \frac{1}{a}\left(\frac{c}{b}+ \mathbb {E}\|X_1^0-X_2^0\|^2e^{-2b t}\right)
	\end{equation}
	for all $t\ge 0$, where $\mathbb {E}$ denotes the \emph{expected value} of the random variable.  The authors additionally prove that this bound can, in general, not be tighter; in particular, one cannot expect $\mathbb {E}(\|X_1(t)-X_2(t)\|^2)\to 0$ as $t\to \infty$.

	The authors also show that if the setup is altered such that $G_1(t,x)\equiv 0$, i.e.~they add noise to a deterministic system and compare it with the noise-free system, then an identical estimate holds with $\frac{c}{b}$ replaced by $\frac{c}{2b}$.
	Further, they extend their theory to interconnected systems of the form
	\begin{equation*}
		d\begin{pmatrix}
			X_1(t) \\
			X_2(t) \\
		\end{pmatrix} = \begin{pmatrix}
			f_1(t,X_1(t),X_2(t)) \\
			f_2(t,X_1(t),X_2(t)) \\
		\end{pmatrix}dt+ \begin{pmatrix}
			G(t,X_1(t)) & 0 \\
			0 & G(t,X_2(t)) \\
		\end{pmatrix}dW(t)
	\end{equation*}
	and give numerous examples of their approach.  However, they do not consider state-dependent metrics $M$.

	In Bouvrie \& Slotine 2019 \cite{BoSl2019contSto}, this last approach is extended under the additional assumptions that $f_1=f_2=f$ are autonomous, i.e.~do not depend explicitly on $t$, that $G_1=G_2=G\colon\R^+_0\times\R^{n\times n}\to \R^n$, and  $G(t,x)G^T(t,x)\ge C_2 I$ for a constant $C_2>0$, i.e.~$GG^T$ is uniformly positive definite. In particular, $W$ is a $2n$-dimensional Wiener process. Then, with $\mu_t$ and $\nu_t$ as the distributions of
	$X_1(t)$ and $X_2(t)$, respectively, one obtains a similar but more informative estimate than \eqref{stocontr1}, namely
	\begin{equation*}
		W_2(\mu_t,\nu_t) \le \frac{1}{\sqrt{a}}\left(\sqrt{\frac{c}{b}}+ W_2(\mu_0,\nu_0)e^{-2b t}\right),
	\end{equation*}
	where
	$$
	W_2(\mu,\nu):= \sqrt{\inf_{\gamma\in \Gamma(\mu,\nu)}\int_{\R^n\times \R^n}\|x-y\|^2_2 \gamma(x,y) dx dy}
	$$
	is the \emph{Wasserstein metric}, measuring the distance between the distributions $\mu$ and $\nu$.   The infimum is taken over all couplings $\Gamma(\mu,\nu)$ of $\mu$ and $\nu$, i.e.~$\gamma\in \Gamma(\mu,\nu)$ is a joint distribution such that
	$$
	\mu(x)=\int_{\R^n}\gamma(x,y)dy\ \ \ \text{and}\ \ \  \nu(y)=\int_{\R^n}\gamma(x,y)dx,
	$$
	and one can think of $W_2(\mu,\nu)$ as the minimal cost of transforming the distribution $\mu$ to the distribution $\nu$.  In this analogy, one can think of the graphs of $\mu$ and $\nu$  picturing piles of earth and the cost of
	moving a unit mass from $x$ to $y$ is $\|x-y\|^2_2$.  The joint distribution $\gamma\in \Gamma(\mu,\nu)$ describes the transport plan and $W_2(\mu,\nu)$ is the square-root of the minimum cost to move earth from the
	pile described by $\mu$ to obtain the pile described by $\nu$.  In computer science, one often refers to the Wasserstein metric, in particular $W_1(\mu,\nu)= \inf_{\gamma\in \Gamma(\mu,\nu)}\int_{\R^n\times \R^n}\|x-y\|_2 \gamma(x,y) dx dy$, as \emph{earth mover's distance}; note that in our analogy $W_1(\mu,\nu)$ is the minimum of the product: amount of earth moved times the (weighted) mean distance of movement.
	
	\subsection{Control systems}
	
	There are numerous publications on the use of contraction analysis in control theory, both regarding \emph{observers}, see e.g.~Aghannan \& Rouchon 2003 \cite{contr2003aghannanrouchon} and Sanfelice \& Praly 2012 \cite{cont2012SP}, and \emph{stabilization} and \emph{control}, see e.g.~Fromion, Monaco \& Normand-Cyrot 1996 \cite{cont1996FMN}, Jouffroy 2003 \cite{control2003Jouffroy}, Slotine 2003  \cite{cont2003Slotine}, Pavlov, van de Wouw \& Nijmeijer 2005/2006 \cite{control2005PWN,control2006PWN}, Pogromski \& Matveev 2016 \cite{control2016PM},
	Wang, Forni, Ortega, Liu \& Su 2017 \cite{control2017WFOL}, Brivadis, Sacchelli, Andrieu, Gauthier \& Serres  2021 \cite{Brivadis2021}, Liu, Xu \& Sun 2021 \cite{control2021LXS}.
	
	Including a detailed discussion of contraction analysis for general control systems in this review would add too much to its length.  We therefore only give a short exposition of \emph{control contraction metrics}  for \emph{control-affine} systems
	\begin{equation}
		\label{affcontsys}
		\dot{x}=\widetilde{f}(t,x,u)= f(t,x)+\underbrace{B(t,x)}_{=(b_1(t,x),\ldots,b_m(t,x))}u=f(t,x)+\sum_{i=1}^m b_i(t,x)u_i
	\end{equation}
	as introduced in Manchester \& Slotine 2017  \cite{control2017MS}.  Assume that
	$u(t,x)=k(t,x)+v(t)$  is smooth
	and denote by $x(t)=\phi(t,t_0,\xi)$ the solution to $\dot{x}=\widetilde{f}(t,x,u)$ with initial value $\xi\in\R^n$ at $t_0$.  The associated (matrix) variational equation is obtained by substituting $x(t)$ for $x$ in \eqref{affcontsys} and taking the derivative with respect to $\xi$
	(recall the matrix-vector product rule $D(A v)=A\cdot Dv + \sum_i v_i Da_i$,  e.g.~\cite[App.~A]{sauer2012})
	\begin{align}
		\frac{\partial}{\partial \xi} \dot{x}(t) &= \frac{d}{dt}\underbrace{\frac{\partial x}{\partial \xi}(t)}_{=:\delta \tilde{x}(t)}=\dot{\delta \tilde{x}}(t)=\frac{\partial }{\partial \xi}\widetilde{f}(t,x(t),\underbrace{u(t,x(t))}_{=k(t,x(t))+v(t)}) \\
		&= \frac{\partial f}{\partial x}(t,x(t)) \delta \tilde{x}(t) +B(t,x(t))\frac{\partial k}{\partial x} (t,x(t))\delta \tilde{x}(t)\n \\
		& \quad +\sum_{i=1}^m u_i(t,x(t)) \frac{\partial b_i}{\partial x}(t,x(t))\delta \tilde{x}(t). \n
	\end{align}
	Note that because we are concerned with the matrix variational equation, here $\delta \tilde{x}(t)\in \R^{n\times n}$.  The connection to $\delta x(t)\in\R^n$ is simple because the  ODE is linear: from the solution $\delta \tilde{x}(t)$ to the matrix variational equation with $\delta \tilde{x}(t_0)=I$, one has $\delta x(t)=\delta \tilde{x}(t)\delta x(t_0)$.
	With
	$$
	A(t,x,u):= \frac{\partial f}{\partial x}(t,x) + \sum_{i=1}^m u_i \frac{\partial b_i}{\partial x}(t,x)\ \ \text{and}\ \ \ K(t,x):=\frac{\partial k}{\partial x} (t,x),%
	$$
	equation \eqref{affcontsys} can be written in the equivalent form
	\begin{equation}
		\dot{\delta x} = \big(A(t,x(t),u(t,x(t)))+B(t,x(t))K(t,x(t))\big)\delta x,
	\end{equation}
	which resembles the classical ODE $\dot{x}=(A+BK)x$  for linear control systems $\dot{x}=Ax+Bu$ with linear feedback control $u=Kx$.
	
	Assume that \eqref{affcontsys} with $u(t,x)=k(t,x)+v(t)$  is contracting in a metric $M(t,x)$ with contraction rate $b=\beta/2>0$.
	With $A:=A(t,x,u(t,x))$, $B:=B(t,x)$, $K:=K(t,x)$, and $M=M(t,x)$ we have by \eqref{contraction0} (see also \eqref{contrrate}) that
	\begin{equation}
		\label{conteqXXX}
		(A+BK)^TM+M(A+BK)+\dot{M} \le - \beta M.
	\end{equation}
	Note that
	\begin{equation}
		\label{conteqXXX2}
		y^TMB=0\ \ \ \text{implies} \ \ \ y^T(A^TM+MA+\dot{M})y \le -\beta y^TMy,
	\end{equation}
	and since $B$ and $M$ are independent of $u(t,x)$, the implication \eqref{conteqXXX2} is a necessary condition for %the existence of the control $u(t,x)=k(t,x)+v(t)$ such that
	\eqref{conteqXXX} to hold true.

	The authors define a uniformly bounded metric $M(t,x)$, i.e.~$\alpha_1 I \le M(t,x) \le \alpha_2 I$ for $0<\alpha_1\le \alpha_2$ to be a \emph{control contraction metric}, if the implication
	\eqref{conteqXXX2} holds true for all $(t,x,u)\in\R\times \R^n\times \R^m$, where $B,K,M$ are as above and $A:=A(t,x,u)$.  Note that the implication must hold true for all $u\in\R^m$, not only for $u=u(t,x)$ for a specific function $u(t,x)$.  They then show, among other things, that if the system
	\eqref{affcontsys} admits a control contraction metric $M(t,x)$, then it is \emph{open-loop controllable} and  \emph{stabilizable} via continuous feedback (universally and exponentially fast).

	% for control systems. More precisely, they showed that the following two properties are equivalent: (i) the input-output map has finite incremental gain and (ii) the unperturbed motion is uniformly asymptotically stable.
	
	%\cite{Brivadis2021}
	%
	%\cite{control2017MS}
	%
	%\cite{control2003Jouffroy} CONTROL?
	%
	%\cite{control2021LXS}
	%
	%\cite{control2017WFOL}
	%
	%\cite{control2016PM}
	%
	%\cite{control2005PWN}
	%\cite{control2006PWN}
	
	%%\cite{cont2003Slotine}
	%many examples and combinations
	%
	%Aghannan and Rouchon 2003  \cite{contr2003aghannanrouchon}.
	%
	%
	%\cite{cont2012SP} input output
	%

	\subsection{Poincar\'e-Bendixson, dimension of attractors, and entropy}
	\label{sec:dimension}

	Smith 1980 \cite{Smith1980} considers a generalization of the Poincar\'e-Bendixson theorem to higher dimensions. He uses a contraction condition with a constant metric, which has two negative and $n-2$ positive eigenvalues, and hence the solutions converge  to a 2-dimensional set, where the Poincar\'e-Bendixson theorem holds. Smith 1987 \cite{Smith1980} continues this work with the following type of assumptions for the autonomous ODE $\dot{x}=f(x)$, where $f$ is locally Lipschitz on its domain $S$: there exists a non-singular matrix $M\in \mathbb R^{n\times n}$ with exactly two negative eigenvalues such that
	$$(x-y)^TM[f(x)-f(y)+\lambda(x-y)]
	\le -\varepsilon \, \|x-y\|^2_2$$
	for all $x,y\in S$
	with constants $\lambda,\varepsilon>0$. The implications are that closed trajectories exist and at least one of them is asymptotically stable.
	
	Sanchez 2009 \cite{Sanchez2009} and 2010 \cite{Sanchez2010} generalizes these concepts in the study of cones of rank 2 and cooperative systems, where a contraction condition of the type
	$$\frac{\partial f^T}{\partial x}(x)M+M\frac{\partial f}{\partial x}(x)+\lambda(x)M<0$$
	holds for a matrix $M\in{\mathcal S}_n$ with 2 negative and $n-2$ positive eigenvalues and $\lambda(x)>0$ is a continuous function, see also the review Burkin 2015 \cite{Burkin2015}.

	Li \& Muldowney 1991 \cite{LM1991} generalize Bendixson's criterion for the nonexistence of closed orbits. They study the eigenvalues of $\frac{\partial f^T}{\partial x}(x) + \frac{\partial f}{\partial x}(x)$ and show that if the sum of the largest two is  positive for all $x$ or if the sum of the smallest two is negative for all $x$, then there are no simple, closed, rectifiable invariant orbits. Similar conditions for the stability and instability of trajectories are derived in Kravchuk, Leonov \& Ponomarenko 1995 \cite{contr1995kravchukleonovponomarenko}.

	\vspace{0.3cm}
	
	The (Hausdorff) \emph{dimension} of an attractor $A$ \cite{dim2011BG} can be bounded by $\dim_H A\le d+s$ if there are an integer $d\in \{0,\ldots,n-1\}$ and $s\in [0,1)$ such that
	\begin{eqnarray}
		\lambda_1(x)+\ldots+\lambda_d(x)+s\lambda_{d+1}(x) < 0\mbox { for all }x\in A, \label{est1}
	\end{eqnarray}
	where $\lambda_1(x)\ge\ldots\ge\lambda_n(x)$ are the eigenvalues of the matrix $\frac{1}{2}(\frac{\partial f^T}{\partial x }(x)+\frac{\partial f}{\partial x }(x))$, %where $Df$ denotes the Jacobian of $f$,
	see Douady \& Oesterl\'e 1980 \cite{dim1980DO}.
	A version only considering contraction in directions $v$, where $v^Tq(x)=0$ and $q^T(x)f(x)\not=0$ is derived in Leonov 1991 \cite{leonov1991}.
	
	This estimate can be improved by means of a Lyapunov-like function, see Leonov \& Boichenko 1992  \cite{dim1992LB} and the book Leonov,  Burkin \&  Shepelyavyi 1996 \cite{contr1996lenonvburkinshepelyavyi}, including applications to, e.g.~the Lorenz and the R\"ossler systems: if  $v\in C^1(\mathbb R^n,\mathbb R)$ is a  function such that (\ref{est1}) is replaced by
	\begin{eqnarray}
		\lambda_1(x)+\ldots+\lambda_d(x)+s\lambda_{d+1}(x)+\dot{v}(x)<0\mbox { for all }x\in A,\label{est2}
	\end{eqnarray}
	%where $v'(x)=\nabla v(x)\cdot f(x)$ denotes the orbital derivative,
	then the same bound on the dimension holds.

	A further generalization in Pogromsky \& Nijmeijer 2000 \cite{dim2000PN} uses a Riemannian metric given by a function   $M\in C^1(\mathbb R^n,{\mathcal S}_n^+)$. The same bound on the dimension holds if the $\lambda_i(x)$ in (\ref{est1}) denote the ordered solutions $\lambda_1(x)\ge\ldots\ge \lambda_n(x)$ to%
	\begin{eqnarray}
		\det\left[\frac{\partial f^T}{\partial x}(x)M(x)+M(x)\frac{\partial f}{\partial x}(x)+\dot{M}(x)-\lambda M(x)\right]=0.\label{est3}
	\end{eqnarray}
	Note that this includes (\ref{est1}) for $M(x) =I$ and (\ref{est2}) for $M(x)=\exp\left(\frac{2v(x)}{d+s}\right)I$ \cite[Cor.~1]{dim2000PN} as well as a constant matrix-valued function $M(x)=M$ as in Leonov 2012 \cite{dim2012Leonov}.
	For similar results on general Riemannian manifolds, see Noack \& Reitmann 1996 \cite{dim1996NR}.
	
	In a similar way, there is a series of results, including  Hartman \& Olech 1962 \cite{contr1962hartmanolech} and Leonov, Burkin \&  Shepelyavyi 1996 \cite{contr1996lenonvburkinshepelyavyi}, to show the following: If
	$$\lambda_1(x)+\lambda_2(x)<0,$$
	with $\lambda$ defined as above in \eqref{est3}, holds for all $x\in \overline{D}$, where $D$ is a bounded, open, simply connected positively invariant set with some additional conditions, then all solutions in $D$ tend to some equilibrium, see Pogromsky \& Nijmeijer 2000 \cite{dim2000PN}. These methods can estimate the dimension without precise localization of the attractor in the phase space, since it is sufficient to apply them to a set containing the attractor, see also the book Boichenko, Leonov \& Reitmann 2005 \cite{contr2005boichenkoleonovreitmann} and Kuznetsov \& Reitmann 2021 \cite{Kuznetsov2021}.
	
	%Finite-time Lyapunov exponents, on the other hand, which are accessible through numerical computation, do not provide a rigorous tool to determine the dimension, see \cite[end of Section 2]{kuznetsov}.
	
	\vspace{0.3cm}
	
	A contraction metric approach to estimate the \emph{topological entropy} for continuous-time systems was introduced in Pogromsky \& Matveev 2010/2011 \cite{POGROMSKY2010,PM2011}, as well as for discrete-time systems in Matveev \& Pogromsky 2016
	\cite{matveev16}. Recently, a converse result for the determination of the so-called \emph{restoration entropy} via a Riemannian metric has been obtained in Kawan, Matveev \& Pogromsky 2020 \cite{kawan-ifac}. The restoration entropy is an upper bound of the topological entropy and has an operational meaning in the context of remote state estimation.

	\subsection{Other extensions}
	\label{sec:other}
	
	Angeli 2002 \cite{contr2002Angeli} presents a framework for incremental stability, extends it to include the \emph{input-to-state stable (ISS)} approach and provides characterizations of incremental stability properties through Lyapunov functions.

	In  Angeli 2009 \cite{contr2009Angeli}, incremental integral input-to-state stability is introduced and a converse result for a continuous Lyapunov function is shown. Furthermore, it is proved that incremental integral ISS implies incremental ISS.
	
	Andrieu, Jayawardhana \& Praly 2016 \cite{man2016AJP} study \emph{transverse stability}, namely a system
	\begin{eqnarray}
		\dot{e}&=&F(e,x)\\\dot{x}&=&G(e,x)
	\end{eqnarray}and the stability of the manifold ${\mathcal E}=\{(e,x)\mid e=0\}$, where $e\in\mathbb R^{n_e}$ and $x\in\mathbb R^{n_x}$.
	This is related to incremental stability by considering a system $\dot{x}=f(x)$ and
	\begin{eqnarray*}
		F(e,x)&=&f(x+e)-f(x)\\
		G(e,x)&=&f(x)
	\end{eqnarray*}
	so that $e$ denotes the error between two solutions $x(t,\xi_1)$ and $x(t,\xi_2)$ with different initial conditions.
	In contrast to the horizontal Finsler-Lyapunov function from \cite{contr2014fornisepulchre}, the contraction condition on the matrix is only required on $\mathcal E$, namely
	\begin{eqnarray}
		\frac{\partial F^T}{\partial e}(0,x)M(x)+M(x)\frac{\partial F}{\partial e}(0,x)+\dot{M}(x)\le -Q\label{ULMTE}
	\end{eqnarray}
	and $p_1I\le M(x)\le p_2 I$ for all $x$; this is called the ULMTE (Uniform Lyapunov Matrix Transversal Equation), where $M(x)\in {\mathcal S}_{n_e}^+$.
	The authors show that the following three properties are equivalent: (i) the manifold $\mathcal E$ is transversally exponentially decreasing, (ii) the transverse linearization along any solution in the manifold is exponentially decreasing, and (iii) there exist positive definite quadratic forms
	whose restrictions to the transverse direction are decreasing. For the case of incremental stability, the first two properties relate to the $e$-(error)-part and the third one is that for every positive definite matrix $Q$ there exists a matrix $M$ as in \eqref{ULMTE} and leads to global results. This is a generalization of a known result for the manifold being an equilibrium. The results are applied to control problems.

	Forni,  Sepulchre \& van der Schaft 2013 \cite{Forni2013} studied \emph{differential passivity} of physical systems using \emph{differential storage functions}, that are also Finsler-Lyapunov functions.

	Incremental stability of \emph{hybrid systems} is considered in Biemond et al.~2018 \cite{biemond2018}.
	
	Generalizations of the contraction analysis have been made to
	\emph{time-periodic} systems in Giesl 2004 \cite{contr2004agiesl} as well as to \emph{almost periodic} systems, see Giesl \& Rasmussen 2008  \cite{contr2008gieslrasmussen}.
	Kawano, Besselink \& Cao 2020
	\cite{time2020KBC} study contraction analysis for \emph{monotone systems}.
	
	Giesl \& Rasmussen 2012  \cite{contr2012gieslrasmussen} studied the non-autonomous ODE $\dot{x}=f(t,x)$ over a \emph{finite time} interval $[0,T]$ and introduced the concept of an \emph{area of (exponential) attraction}, which is a set $G\subset [0,T]\times \mathbb R^n$ such that all solutions in $G$ are (exponentially) attractive. In contrast to the domain of attraction, which depends on one particular solution, the area of attraction does not depend on one special solution. They  showed that $G$ is an area of exponential attraction if and only if  there exists a Riemannian metric $M\colon [0,T]\times \mathbb R^n\to {\mathcal S}_n^+$ which satisfies
	\begin{align*}&
		\frac{1}{2}	\max_{w\in\mathbb R^n,w^TM(t,x)w=1}
		w^T\left(\frac{\partial f^T}{\partial x}(t,x)M(t,x)+M(t,x)\frac{\partial f}{\partial x}(t,x)+\dot{M}(t,x)\right)w \\
	&\le  -\nu<0
	\end{align*}
	for all $(t,x)\in G$, where $G\subset [0,T]\times \mathbb R^n$ is connected and invariant,
	together with $M(0,x)=M(T,x)=I$ for all $x\in \mathbb R^n$.
	The rate of exponential attraction is  linked to $\nu$.

	\vspace{0.3cm}

	Ngoc \& Trinh 2018 \cite{del2018NT} derive contraction criteria for \emph{delay differential equations} (functional differential equations) of the form
	$$\dot{x}(t)=f(t,x(t),x_t),\quad t\ge \sigma,$$
	where $x_t(\cdot)$ is defined by $x_t(\theta):=x(t+\theta)$, $\theta\in [-h,0]$. Hence, $\dot{x}$ at time $t$ depends on the value of the solution $x$ in $[t-h,t]$.
	A criterion for the exponential contraction between two solutions is given, which is based on a uniform bound on
	$\frac{\partial f}{\partial x}(t,x,\varphi)$, independent of $x$ and $\varphi$.
	Furthermore, if $f$ is periodic in the first argument, a similar criterion is presented that shows the existence of a unique periodic solution, which is globally exponentially stable.
	
	In Ngoc, Trinh, Hieu \& Huy 2019 \cite{del2018NTLH}, criteria for the contraction of difference equations with time-varying delays are presented.

	\vspace{0.3cm}

	The book Leine \& van de Wouw 2008 \cite{Leine2008} studies contraction analysis for \emph{non-smooth systems}.
	It contains a literature survey of stability and convergent systems in non-smooth systems (Chapter 1.5) and considers convergent systems in Chapter 8.
	To study non-smooth systems, the book uses the framework of measure differential inclusions, which includes systems with state discontinuities.
	The authors introduce the concepts of stability, incremental stability and convergence for these systems. Convergent systems are studied in Chapter 8, based on Leine \& van den Wouw 2008 \cite{leine2008b}, where the authors consider maximal monotone systems and use, e.g., incremental stability arguments to show convergence (Theorem 8.7).
	An early reference for non-smooth systems is Reissig 1954 \cite{Reissig1954}, deriving conditions for the convergence of all solutions to a unique solution.

	Giesl 2005 \cite{contr2005giesl}  considers a one-dimensional time-periodic and nonsmooth ODE $\dot{x}=f(t,x)$ with $x\in \mathbb R$, $f(t+T,x)=f(t,x)$ and $f$ is $C^1$ outside of the line $x=0$. He derives a contraction condition for the existence of an exponentially stable periodic orbit based on the distance between $(t,x)$ and $(t,y)$ given by
	$$A(t)=e^{W(t,x)}|y-x|,$$
	where $W(t,x)$ is a time-periodic function, which is smooth outside of the line $x=0$. The contraction conditions are
	$$\frac{\partial f}{\partial x}(t,x)+\dot{W}(t,x)\le-\nu<0$$
	for $x\not =0$ together with two jump conditions for points $(t,0)$ where solutions pass the line $x=0$.
	In Giesl 2007 \cite{giesl-nonsmooth-nec}, a converse theorem is proved, showing the existence of such a metric. Stiefenhofer \& Giesl 2019 \cite{SG2019,SG2019-2} generalize the sufficient condition to two spatial dimensions.

	Dahlquist 1958 \cite{num1958Dahlquist}
	applies contractivity in numerics for solving ODEs.  Desoer 1972 \cite{num1972DH} uses that contractive systems also satisfy an ISS property to derive bounds between trajectories of a continuous-time (contractive) system and its discretization. This is applied to the numerical solution of the ODE system.
	In Wensing \& Slotine 2020 \cite{WeSL2020contrGrad}, the convergence of the \emph{gradient descent} algorithm for minimization is studied using contraction analysis and geodesic convexity, and classical results based on convexity of the objective function are generalized.
	Forni \& Sepulchre 2019 \cite{FoSe2019domana} develop differential dissipativity theory for the \emph{dominance analysis} of nonlinear systems.
	
	\section{Numerical construction methods}
	\label{sec:num}
	
	In this last section, we discuss construction methods for contraction metrics. We classify them by the numerical methods employed and briefly describe them.
	Except for the \emph{subgradient method}, which uses optimization on a matrix manifold, these methods have counterparts in the numerical computation of Lyapunov functions, see e.g.~Giesl \& Hafstein 2015 \cite{GiHa2015review},
	but are more involved because a matrix-valued contraction metric, rather than a real-valued Lyapunov function, is computed.

	%\subsection{Analytical construction}
	
	\subsection{Collocation}
	\label{sec:collo}
	
	In this method, a specific contraction metric is characterized by a partial differential equation (PDE). Then meshfree collocation is used to approximately solve this PDE in a given compact set. The approximation, if sufficiently close, is itself a contraction metric
	because it fulfills the required inequalities.
	
	In this section,    $M(x)$ stands for a specific contraction metric, satisfying the equation \eqref{eq1}, while $S(x)$ denotes its approximation, which fulfills the inequalities \eqref{S2}, \eqref{S1}.
	For $S(x)$ to be a contraction metric of the autonomous system $\dot{x}=f(x)$, we require
	\begin{eqnarray}
		\frac{\partial f^T}{\partial x}(x)S(x)+S(x)\frac{\partial f}{\partial x}(x)
		+\dot{S}(x)&<&0\label{S2}\\
		S(x)&>&0\label{S1}
	\end{eqnarray}
	for all $x\in K$, where $K\subset \mathbb R^n$ is a compact set.
	Converse theorems show the existence of a contraction metric $M(x)$ satisfying
	\begin{eqnarray}
		\frac{\partial f^T}{\partial x}(x)M(x)+M(x)\frac{\partial f}{\partial x}(x)
		+\dot{M}(x)&=&-C\label{eq1}\\
		M(x)&>&\varepsilon I,\label{eq2}
	\end{eqnarray}
	see Giesl 2015 \cite[Th.~4.4]{contr2015giesl}.
	Note that \eqref{eq1} with a given matrix $C\in {\mathcal S}_n^+$ has a unique solution which satisfies \eqref{eq2} in every compact set within the basin of attraction for some $\varepsilon >0$. Hence, a sufficiently close  approximation $S$ to $M$ will satisfy \eqref{S2} and \eqref{S1} as well.
	
	Let us introduce meshfree collocation in more detail. For a general overview, see Buhmann 2003 \cite{rbf2003buhmann} or Wendland 2005 \cite{rbf2005wendland}; for the case of matrix-valued functions, see Giesl \& Wendland 2018
	\cite{contr2018GW}.
	
	In general, one considers a linear operator $L\colon H\to H$ acting on a Hilbert space $H$. Here, $H=H^\sigma(\Omega,{\mathcal S}_n)$ is a \emph{reproducing kernel Hilbert space} (RKHS), which  consists of the same functions as the Sobolev space of matrix-valued functions, and their norms are equivalent. The definition of the kernel $\Phi:\Omega\times\Omega\to \cL(\cS_n)$, where $\cL$ denotes the linear space of
	linear and bounded operators $L:\cS_n\to \cS_n$, is based on the kernel for the RKHS $H^\sigma(\Omega,\mathbb R)$ of scalar-valued functions.
We consider the linear operator
	$$LM(x)=\frac{\partial f^T}{\partial x}(x)M(x)+M(x)\frac{\partial f}{\partial x}(x)+\dot{M}(x).$$ The PDE  to be solved can be written as
	$$LM(x)=r(x),$$
	where $r\colon\R^n\to \cS^+_n$ is a given function.
	We choose a finite set of collocation points $x_1,\ldots,x_N\in \Omega$ and define the corresponding linear operators $\lambda_k^{(i,j)}\in H^*$ by $\lambda_k^{(i,j)}(M)=(LM(x_k))_{ij}$, where $k=1,\ldots,N$ and $1\le i\le j\le n$.
	
	The approximation can now be formulated as a generalized interpolation problem, namely finding the norm-minimal function which solves the PDE at the collocation points, i.e.
	$$\min_{S\in H}\{\|S\|_H\mid  \lambda_k^{(i,j)}(S)=(r(x_k)_{ij}) \text{ for }k=1,\ldots,N \text{ and }1\le i\le j\le n\}.$$
	
	The solution of the minimization problem is unique and given by a linear combination of the Riesz representers  of the $\lambda_k^{(i,j)}$,
	where the coefficients are chosen such that the interpolation conditions hold.  Note that the Riesz representer of $\lambda\in H(\Omega,\cS_n)^*$ in the RKHS is given by
	\begin{eqnarray}
		v_\lambda(x) = \sum_{1\le \mu\le\nu\le n}
		\lambda(\Phi(\cdot,x)E^s_{\mu\nu})E^s_{\mu\nu}, \qquad x\in\Omega,
		\label{Rieszr}
	\end{eqnarray}
	where $(E^s_{\mu\nu})_{1\le \mu\le\nu\le n}$ is an orthonormal basis of $\cS_n$ and $\Phi$ is the kernel.
	The interpolation condition is a system of linear equations with collocation matrix
	of size $Nn(n+1)/2$. To find the coefficients, and thus the solution, this system of linear equations needs   to be solved.
	%
	%\times N}$, $a_{ij}=\lambda_i(v_j).$%_{i,j=1,\ldots,N}$ and $$a_{ij}=\lambda_i(v_j).$$
	%If $H$ is a reproducing kernel Hilbert space with kernel $\Phi\colon \Omega\times \Omega\to H$, then the Riesz representers are given by $v_j=\lambda_j^y\Phi(x,y)$ and we have $a_{ij}=\lambda_i^x\lambda_j^y\Phi(x,y)$ and $A\in\cS^+_N$. %the collocation matrix becomes the symmetric matrix
	%%$$a_{ij}=\lambda_i^x\lambda_j^y\Phi(x,y).$$
	
	The computation of a contraction metric for an {\it equilibrium} using meshfree collocation was studied in Giesl \& Wendland 2019 \cite{contr2019GW}. The construction was achieved by approximating the contraction metric $M$ satisfying \eqref{eq1}. If $\Omega $ is in the basin of attraction of the equilibrium and $f$ is sufficiently smooth, then the approximation $S$ satisfies
	\begin{eqnarray}
		\sup_{x \in K} \|M(x)-S(x)\|_2 \le
		C_1 \|L(M)-L(S)\|_{L^\infty(\Omega,\cS_n)} \le C_2
		h_{X,\Omega}^{\sigma-1-n/2}
		\|M\|_{H^\sigma(\Omega,\cS_n)}
		\label{RBF-est}
	\end{eqnarray} on every positively invariant compact set $K\subset \Omega$ with some constants $C_1,C_2>0$. Here, $h_{X,\Omega}:=\min_{x\in X}\sup_{y\in \Omega} \|y-x\|_2$ is the so-called fill distance of the set of collocation points $X=\{x_1,\ldots,x_N\}$ in $\Omega$ \cite{contr2018GW}.
	It follows that $S$ fulfills \eqref{S2} and \eqref{S1} and is a contraction metric for the system if $h_{X,\Omega}$ is small enough, i.e.~if the collocation points are dense enough.

	%For  $ W \in  C^k(\Omega; \R^{n \times n}) $, where $\cD\subset \R^n$ is a non-empty open set and $ \cR $ is  $ \R, \R^n, \Sb^{n\times n}, $ or $\R^{n \times n}$, we define the $C^k$-norm  as
	%\begin{equation}
	%\norm{W}{C^k(\cD;\cR) } := \sum_{| {\balpha} | \leq k} \, \sup_{\bx \in \cD} \norm{\partial^{\balpha} W(\bx)}{2},
	%\end{equation}
	%where $ \balpha\in\N_0^n $ is a multi-index and $ |\balpha| := \sum_{i=1}^{n} \alpha_i$.  When all  $D^{\balpha} W$ can be continuously extended to $\overline{\cD}$ for all $|\balpha|\le k$, the $C^k$-norm is also defined on $\overline{\cD}$ with the same formula.

	A contraction metric for systems with an exponentially stable {\it periodic orbit} was first computed for two-dimensional systems, $\dot{x}=f(x)$, $f\colon \mathbb R^2\to\mathbb R^2$, where a metric of the form  $M(x)=\exp(2W(x))I$ exists and thus a scalar-valued function $W\colon \mathbb R\to\mathbb R$ needs to be found. In Giesl 2007 \cite{rbf2007agiesl}, the periodic orbit and the Floquet exponent were approximated numerically and then $W$ was found by solving a first-order PDE for $W$ using meshfree collocation.
	In Giesl \& McMichen 2016 \cite{mcmichen}, no information about the periodic orbit is required and $W$ is found by solving a second-order PDE using meshfree collocation.
	
	In higher dimensions, the metric can in general not be characterized by a scalar-valued function.
	Giesl 2009
	\cite{rbf2009giesl} computes the local metric by numerically approximating the periodic orbit and computing a metric $M(x)$ using the first variational equation.
	
	A contraction metric was characterized as the solution of a  matrix-valued PDE, see
	Giesl 2021 \cite{contGiesl2021} and \eqref{peri1}, \eqref{peri2}. The computation using meshfree collocation was achieved in  Giesl 2019	
	\cite{contr2019aGiesl}.
	
	Finally, Giesl \& McMichen 2018 \cite{fin2018GM} consider the one-dimensional, non-autonomous system $\dot{x}=f(t,x)$ on a {\it finite time interval}, where $x\in \mathbb R$ and $t\in [T_1,T_2]$. The area of exponential attraction, see Section \ref{sec:other} and \cite{contr2012gieslrasmussen}, can be determined by a contraction metric. They show that the contraction metric can be chosen of the form $M(t,x)=\exp(2W(t,x))$, where the scalar-valued function $W$ satisfies a second-order PDE involving the second orbital derivative $\ddot{W}(t,x)$, i.e.~the orbital derivative of the orbital derivative $\dot{W}(t,x)$,
	$$
	\ddot{W}(t,x)=-\left(\frac{\partial^2f}{\partial x \partial t}(t,x)+\frac{\partial^2f}{\partial x^2 }(t,x)f(t,x)\right)
	$$ with boundary values $W(T_1,x)=W(T_2,x)=0$. Meshfree collocation is used to solve this equation numerically and to compute a contraction metric.

	\subsection{Linear matrix inequalities}
	
	\emph{Linear matrix inequalities} (LMIs) are closely linked to \emph{semidefinite optimization problems} of the form:  \\
	
	Given $C,B_i\in\cS_n$, $i=1,\ldots, N,$ and a vector $c\in\R^N$, determine a vector $q\in\R^N$ to
	\begin{equation}
		\label{LMIform}
		\text{minimize}\ \ \ c^Tq \ \ \ \text{subject to}\ \ \ \sum_{i=1}^Nq_iB_i+C\ge 0,
	\end{equation}
	i.e.~the linear functional $q \mapsto c^Tq $ is minimized under the given conditions on the vector $q$.
	If $c=0$, the problem \eqref{LMIform} is called a \emph{feasibility problem} and any $q\in\R^N$ such that the constraints are fulfilled, if there are any, is a solution to the problem.
	In a sense, LMI problems are a generalization of the Lyapunov equation $A^T M+MA=-C$, for given $A\in\R^{n\times n}$, $C\in \cS^+_n$.
	%, which is a special case of the Sylvester equation.
	Indeed, LMIs are in much use in control theory and the literature on the use of LMIs is extensive, see e.g.~the monographs  \cite{BGFB1994LMIincontrol,sosphdthesis2000parrilo,chesi2010LMIsurvey,chesibook} for an introduction.
	Note that an LMI problem like `find a symmetric matrix fulfilling $M>0$ and $A^TM+MA<0$' can be written in the form \eqref{LMIform} and % and we will use LMI problem and the semidefinite programming problem
	LMIs can be routinely solved using \emph{semidefinite programming}, which is a subclass of \emph{convex optimization} \cite{BoVa2004conopt}. There are numerous computational packages available for different platforms, e.g.~MATLAB's Robust Control Toolbox (formerly LMI toolbox), CVX, YALMIP, and Mosek \cite{S98guide,cvx}.
	
	Related to LMIs are so-called \emph{bilinear matrix inequalities} (BMIs), where additionally  matrices $A_{ij}\in\cS_n$, $i,j=1,\ldots,N$, are given and the objective is to determine a vector $q\in\R^N$ solving the problem
	\begin{equation}
		\label{BMIform}
		\text{minimize}\ \ \ c^Tq \ \ \ \text{subject to}\ \ \ \sum_{i,j=1}^Nq_iq_jA_{ij}+\sum_{i=1}^Nq_iB_i+C\ge 0.
	\end{equation}
	Such problems are in general not convex and are much harder to solve \cite{vABr2000lmibmi}. In addition to the computation of contraction metrics, BMIs have also been considered for the computation of Lyapunov functions for SDEs \cite{Ha2019BMI}.

	In the following sections, we will list several methods developed to reformulate the computation of a contraction metric to LMIs and one example of the use of BMIs.

	\subsubsection{Sum of squared polynomials}
	
	SOS is originally an approach to show that a polynomial is positive, see Parrilo 2000 \cite{sosphdthesis2000parrilo}.
	A multivariate polynomial $p(x)\in\mathbb R[x]$, where $\mathbb R[x]$ denotes the polynomials in $x\in\mathbb R^n$, is SOS (sum of squares) if it can be written as
	$$p(x)=\sum_{i=1}^m p_i^2(x),$$
	where $p_1(x),\ldots,p_m(x)\in \mathbb R[x]$. It is clear that an SOS polynomial is non-negative, but not every non-negative polynomial is SOS, see e.g.~\cite{motzkin} for a counterexample.
	
	A polynomial is SOS if and only if there exists a positive semidefinite matrix $Q$ such that
	\begin{eqnarray}p(x)&=&Z^T(x)QZ(x),\label{SOS-link}
	\end{eqnarray}
	where $Z(x)$ is a vector of monomials of degree less than or equal to half the degree of $p$.
	The link to LMIs of the form \eqref{LMIform} is to write $Q=\sum_iq_iB_i$ where $B_i$ is a basis of $\cS_n$.
	
	Let us explain \eqref{SOS-link} with an example: consider a polynomial $p(x_1,x_2)$ of degree 4.  Then $p(x_1,x_2)$ is an SOS polynomial if there is a $Q\ge 0$ such that
	$$
	p(x_1,x_2) = Z^TQZ= \begin{pmatrix}
		1 & x_1 & x_2 & x_1x_2 & x_1^2 & x_2^2 \\
	\end{pmatrix}Q
	\begin{pmatrix}
		1 \\ x_1 \\ x_2 \\ x_1x_2 \\ x_1^2 \\ x_2^2 \\
	\end{pmatrix}.
	$$
	Since $Q\ge 0 $, it can be written as $Q=O^TDO$, where $D=\text{diag}(d_1,d_2,\ldots,d_6)$ is a diagonal matrix with non-negative entries and $O$ is an orthogonal matrix.
	The entries $p_i=p_i(x_1,x_2)$ of the vector $OZ=(p_1,p_2,\ldots,p_6)^T$ are polynomials in $(x_1,x_2)$ and
	$$
	p(x_1,x_2)=Z^TO^TDOZ = [OZ]^TD[OZ] =\sum_{i=1}^6 \left(\sqrt{d_i}\,p_i(x_1,x_2)\right)^2.
	$$
	Note that these conditions hold globally in $\mathbb R^n$.
		It is also possible to search for SOS polynomials on compact domains by use of the
		Positivstellensatz, i.e.~an algebraic description of all polynomials that are positive, see also Hilbert's 17th problem.

	This has been used, e.g., for the construction of Lyapunov functions by solving LMI problems numerically, see Chesi 2010 \cite{chesi2010LMIsurvey} and Anderson \&  Papachristodoulou 2015 \cite{AnPa2015} for an overview. There is also a free toolbox for MATLAB available to compute SOS Lyapunov functions
	for polynomial systems, see  Papachristodoulou et al. 2013 \cite{PAVPSP2013sostools}.
	
	For the construction of SOS contraction metrics, see  Aylward, Parrilo \& Slotine 2008 \cite{contr2008aylwardparriloslotine}, an SOS matrix  \cite{GP2004} needs to be defined: a symmetric matrix with polynomial entries $S(x)\in \mathbb R[x]^{n\times n}$ is an SOS matrix if the scalar polynomial
	$$y^TS(x)y$$ is SOS in $\mathbb R[x,y]$, where $y\in\mathbb R^n$.
	% Similarly to the scalar case, a matrix is SOS if and only if there is a matrix $T(x)\in\mathbb R[x]^{p\times m}$ such that
	%$$S(x)=T(x)^T T(x).$$
	To deal with strict inequalities, they define $S(x)$ to be a strict SOS matrix if $$S(x)-\varepsilon I$$ is an SOS matrix for some $\varepsilon>0$.
	
	The idea is now to find a matrix-valued function $M(x)\in{\mathcal S}_n$ such that the entries are polynomials of a fixed maximal degree $d$ and $M(x)$ as well as $$-R(x):=-\left(\frac{\partial f^T}{\partial x}(x)M(x)+M(x)\frac{\partial f}{\partial x}(x)+\dot{M}(x)\right)$$ are strict SOS matrices.
	Alternatively, to show contraction with a contraction rate $\beta/2$, one requires that  $$-\left(\frac{\partial f^T}{\partial x}(x)M(x)+M(x)\frac{\partial f}{\partial x}(x)+\beta M(x)+\dot{M}(x)\right)$$ and $M(x)$ are strict SOS matrices.
	In more detail, one writes $M(x)$ with unknown coefficients (polynomials) and derives the conditions in the unknown coefficients. The constraints translate into SOS constraints on the scalar polynomials $y^TM(x)y $ and $-y^TR(x)y$. If the SOS solver finds a solution, then
	it can be used to generate
	%we have found
	a contraction metric.
	
	Aylward, Parrilo \& Slotine 2008 \cite{contr2008aylwardparriloslotine} apply the method to an example and show that the degree of the polynomials needs to be sufficiently large to construct a contraction metric; the larger the degree, the better the bound on the rate of convergence. Furthermore, they show that the constructed contraction metric is a valid contraction metric for a perturbed system and they analyze the range of perturbations.
	
	SOS was also used to construct a contraction metric for a periodic orbit of an autonomous system in Manchester  \& Slotine 2014  \cite{contr2014manchesterslotine}. To remove the condition $v^TM(x)f(x)=0$ from  \eqref{cont-per1} (in autonomous form), they define $W(x)=M^{-1}(x)$ and $w=M(x)v$, so that \eqref{cont-per1} becomes
	\begin{eqnarray*}
		w^T\left(W(x)\frac{\partial f^T}{\partial x}(x)+\frac{\partial f}{\partial x}(x)W(x)+W(x)\dot{M}(x)W(x)\right)w
		\le -\beta w^TW(x){w}.
	\end{eqnarray*}
	Since the condition
	$0=v^TM(x)f(x)=w^Tf(x)$ is equivalent to $w^Tf(x)f^T(x)w\le 0$ and the overall condition can be written as
	\begin{eqnarray*}
		H(x):=
		W(x)\frac{\partial f^T}{\partial x}(x)+\frac{\partial f}{\partial x}(x)W(x){-\dot{W}(x)}+\beta W(x)-\rho(x) f(x)f^T(x)
		\le 0
	\end{eqnarray*}
	with a function $\rho(x)\ge 0$.
	They now apply SOS to construct a contraction metric, namely to find a polynomial matrix-valued function $W(x)$ and a scalar-valued function $\rho(x)$ such that $W(x)$, $H(x)$ and $\rho(x)$ are SOS (on a compact set, using the Positivstellensatz).
	
	In August \& Barahona 2011 \cite{august2011}, the authors present sufficient conditions for global complete synchronization of coupled identical oscillators, which are linear matrix inequalities, based on contraction analysis; the point-dependent metric is then calculated using SOS.

	\subsubsection{Reaction-diffusion equation}
	Arcak 2011 \cite{graph2011Arcak} studied the reaction-diffusion equation
	\begin{equation}
		\label{diffusionEQ}
		\dot{x}=f(x)+D(-\Delta) x,\ \ \ \text{where} \ \ D\in\R^{n\times n}\ \ \text{and}\ \ (-\Delta) x=-\left(\Delta x_1,\ldots,\Delta x_n\right)^T
	\end{equation}
	is the vector Laplacian,
	subject to  Neumann boundary conditions $\nabla x_i(t,\xi)\cdot n(\xi)=0$, $\xi\in \partial \Omega$, on a bounded domain $\Omega$.
	Note that since we write the Laplacian in its positive definite form $-\Delta$, our $D$ is $-D$ in \cite{graph2011Arcak}.
	With $\lambda_2$ as the second smallest eigenvalue of the operator $-\Delta$ on $\Omega$ with vanishing Neumann condition, it is shown
	that if there exists a convex set $\cX\subset \R^n$, $M\in \cS^+_n$, and $\varepsilon>0$, such that
	\begin{align}
		&D^TM+MD\le 0\ \ \text{and}\n\\
		\label{diffusLMI1} &\left(\frac{\partial f}{\partial x}(x)+\lambda_2D\right)^TM+M\left(\frac{\partial f}{\partial x}(x)+\lambda_2D\right) \le -\varepsilon I \ \ \text{for all $x\in\cX$},
	\end{align}
	then every classical solution $x\colon\R_0^+ \times \Omega \to \cX$ to \eqref{diffusionEQ} in $\cX$ converges exponentially fast in time $t$ to its average in the $L_2(\Omega)$-norm.
	In formulas, there exist constants $b,C>0$ such that
	$$
	%\left(\int_\Omega (x(t,\xi)-\overline{x(t)})^2 d^n\xi\right)^{\frac12} \le Ce^{-bt},\ \ \text{where}\ \ \overline{x(t)}:=\frac{1}{|\Omega|}\int_\Omega x(t,\xi) d^n\xi.
	\int_\Omega (x(t,\xi)-\overline{x}(t))^2 d^n\xi \le Ce^{-bt},\ \ \text{where}\ \ \overline{x}(t):=\frac{1}{|\Omega|}\int_\Omega x(t,\xi) d^n\xi.
	$$
	Note that the smallest eigenvalue $\lambda_1$ of $-\Delta$ is zero and the second smallest eigenvalue $\lambda_2$ is related to the connectedness of the domain $\Omega$ and is maximized for a ball-shaped $\Omega$.

	The LMI  conditions  \eqref{diffusLMI1} are then made tractable
	for a numerical procedure by assuming
	$$
	\frac{\partial f}{\partial x}(x)\in  \underbrace{\conv\{Z_1,\ldots,Z_q\}}_{:=\{\sum_{i=1}^q\mu_iZ_i,\ \mu_i\ge 0,\ \sum_{i=1}^q\mu_i=1\}}   + \underbrace{\operatorname{cone}\{S_1,\ldots,S_m\}}_{:=\{\sum_{i=1}^m\mu_iS_i,\ \mu_i\ge 0\}}
	\ \ \ \text{for all $x\in \cX$}.
	$$
	Then condition \eqref{diffusLMI1} is fulfilled for some $\varepsilon>0$ if
	\begin{align*}
		(Z_i+\lambda_2D)^TM+ M(Z_i+\lambda_2D) <0 \ \ \text{and}\ \ S_j^TM+MS_j\le 0
	\end{align*}
	for $i=1,\ldots,q$ and $j=1,\ldots,m$.

	Alternatively, one can assume that
	$$
	\frac{\partial f}{\partial x}(x)\in  \underbrace{\operatorname{box}\{A_0,b_1c_1^T,\ldots,b_\ell c_\ell^T\}}_{:=\{A_0+\sum_{i=1}^\ell \mu_i b_i c_i^T ,\ 0 \le \mu_i\le 1\}},\ \ A_0\in\R^{n\times n},\ b_i,c_i\in\R^n, \ \ \ \text{for all $x\in \cX$}.
	$$
	Then, condition \eqref{diffusLMI1} is fulfilled for some $\varepsilon>0$ if there is an
	$$
	\cP=\begin{pmatrix}
		M &  &  &  \\
		& q_1 &  &  \\
		&  & \ddots &  \\
		&  &  & q_\ell\\
	\end{pmatrix}\in \R^{(n+\ell)\times(n+\ell)} ,\ \ \ q_i>0\ \  \text{for $i=1,\ldots,\ell$,}
	$$
	such that with $B=(b_1,\ldots,b_\ell)\in \R^{n\times \ell}$ and $C=(c_1,\ldots,c_\ell)\in \R^{n\times \ell}$, i.e.~the $b_i$ and $c_i$ are the columns of the matrices $B$ and $C$, we have
	$$
	\cP\begin{pmatrix}
		A_0+\lambda_2D & B \\
		C^T & -I \\
	\end{pmatrix}+\begin{pmatrix}
		A_0+\lambda_2D & B \\
		C^T & -I \\
	\end{pmatrix}^T\cP<0.
	$$
	Three examples of the use of these LMIs are presented using analytical and numerical solutions.  Further, an example is presented where the conditions  $D^TM+MD\le 0$ and
	\eqref{diffusLMI1} cannot be simultaneously fulfilled.  Finally, the theory is expanded to compartmental ODEs of the form
	$$
	\dot{x}^k=f(x^k)+D\sum_{j\in\cN_k}(x^j-x^k),\ \ D\in\R^{n\times n},\ x^k\in\R^n,
	$$  where $\cN_k\subset \{1,\ldots,N\}$
	is the set of nodes adjacent to node $k$. % The interpretation is that we have $N$ identical systems $\dot{x}=f(x)$ (nodes), where some are interconnected through diffusion terms (the adjacent nodes).
	The term $\sum_{j\in\cN_k}(x^j-x^k)$ corresponds to a graph Laplacian and the theory can be adapted to this case to show an exponential rate of synchronization of $x^k(t)$ and $x^j(t)$ as $t\to \infty$ for all $j,k\in\{1,\ldots,N\}$.  For a more general approach to synchronization using contraction, see Russo, di Bernardo \& Sontag 2013 \cite{graph2013RBS}.
	
	\subsubsection{Switched systems}
	
	In Pavlov, Pogromsky, van de Wouw \& Nijmeijer 2007 \cite{PPWN2007}, the convergence of piecewise affine systems with bounded, piecewise continuous inputs is studied and sufficient and necessary LMI conditions are derived.
	The  approach is motivated by the system $\dot x =Ax+w(t)$ with a Hurwitz $A\in\R^{n\times n}$ and bounded, piecewise continuous $w\colon\R\to \R^m$. Its solution with $x(0)=x_0$ is $x(t)=e^{At}(x_0+\int_0^te^{-A\tau}w(\tau)d\tau)$ and one easily verifies that the only solution bounded on $\R$ is the one with initial state $x_0=\int_{-\infty}^0 e^{-A\tau}w(\tau)d\tau$.  Further, by considering the Jordan normal form of  $A$, one  has for $t\ge0$ and any solution $x^*(t)$ with $x^*(0)=x^*_0$ that
	$$
	\|x^*(t)-x(t)\|_2 \le \|e^{At}(x^*_0-x_0)\|_2 \le Ce^{-bt}\|x^*_0-x_0\|_2,
	$$
	where $-b<0$ is larger than the real part of any eigenvalue of $A$, which are all negative since $A$ is Hurwitz. Thus, for every admissible input $w$, there is a unique bounded solution $\bar x(t)$ to $\dot x =Ax+w(t)$ and this solution contracts all other solutions exponentially fast as $t\to \infty$.
	%Here one can show that if the existence of a unique bounded solution $\bar x(t)$ defined on the whole of $\R$, that exponentially extracts all other solutions.
	A system with these properties is called \emph{exponentially convergent}; see also Section \ref{sec:difftypesstab}.
	Alternatively, one can consider this unique bounded solution as a non-autonomous \emph{pullback attractor}, see Kloeden \& Rasmussen \cite{nonaut2011kloedenrasmussen}.
	
	With matrices $A_i\in\R^{n\times n}$, vectors $b_i\in\R^n$, and input $w\in\R^n$, the authors study piecewise affine systems of the form
	\begin{equation}
		\label{PavSys}
		\dot x =A_ix+b_i+w, \ \ x\in \Lambda_i,
	\end{equation}
	with $\R^n$ divided into a finite number of polyhedral cells $\Lambda_i$ by a finite number of hyperplanes $H_j:=\{x\in\R^n\mid v_j^Tx+h_j=0\}$, $v_j\in\R^n$ and $h_j\in\R$.
	
	The authors use Lyapunov functions $V(x_1,x_2)=(x_1-x_2)^TP(x_1-x_2)$, $P\in \cS^+_n$, of two arguments to study stability.  Note that this Lyapunov function measures the distance between two solutions $x_1(t)$ and $x_2(t)$.
	They first show that the LMIs
	\begin{equation}
		\label{paarp}
		A_i^TP+PA_i<0 \ \ \ \text{for all $i$}
	\end{equation}
	are sufficient for the exponential convergence of the system \eqref{PavSys}, if its right-hand side is continuous.  Further, they show that the right-hand side is continuous, if and only if whenever two cells
	$\Lambda_i$ and $\Lambda_j$ intersect in a hyperplane $H_k$, one has
	$$
	h_k(A_i-A_j)=(b_i-b_j)v_k^T.
	$$
	If the system \eqref{PavSys} does not have a continuous right-hand side, then \eqref{paarp} is not a sufficient condition and, in the simplified case of a bimodal system with one hyperplane $v^Tx=0$ and $i=1$ for $v^Tx\ge 0$ and $i=2$ otherwise,
	they show that the following LMI  is sufficient for $P\in \cS^+_n$ to define a Lyapunov function of two arguments as above, that establishes exponential convergence:
	Define $\Delta A:=A_1-A_2$ and $\Delta b:=b_1-b_2$.  The (constrained) variables of the  LMI are $P\in \cS^+_n$ and real numbers $\beta>0$, $\gamma \ge 0$. The additional constraints are
	$$
	\begin{pmatrix}
		PA_1+A_1^TP+\beta I & P\Delta A-\frac{1}{2}vv^T \\
		&\\
		\Delta A^TP-\frac{1}{2}vv^T & -vv^T \\
	\end{pmatrix} \le 0\ \ \ \text{and $P\Delta b=-\gamma v$.}
	$$
	Although these conditions can obviously be implemented in software, the authors do not provide a numerical example.

	In Fiore, Hogan \& Bernardo 2016 \cite{swi2016FHB}, these results were advanced using regularization and numerous analytically  solved examples were presented.

	\subsubsection{Differential-algebraic systems}
	In Nguyen,  Vu, Slotine \& Turitsyn 2021 \cite{NVST2021contrDAE}, the contraction in nonlinear \emph{differential-algebraic equations} (DAEs)  was studied.  DAEs  are concerned with system descriptions of the form
	\begin{equation}
		\label{DAEgeneral}
		\dot{x}=f(x,y),\ \ \ g(x,y)=0,
	\end{equation}
	where $x\in \R^n$, $y\in \R^m$, and $g:\R^{n+m}\to \R^m$ models dynamics that are assumed to be infinitely fast in comparison to the dynamics modelled by $f:\R^{n+m}\to\R^n$.
	%They study contraction in for general $p$-norms
	Given a point $(x_0,y_0)\in\R^{n+m}$   such that $g(x_0,y_0)=0$ and $\frac{\partial g}{\partial y}(x_0,y_0)\neq 0$, the implicit function theorem ensures the existence of a unique function $\phi$ such that $\phi(x_0)=y_0$ and $g(x,\phi(x))=0$ in a neighbourhood of $(x_0,y_0)$.  Inductively applying this argument, one obtains a maximal
	domain
	$$
	\cR=\left\{x\in\R^n\mid g(x,\phi(x))=0\ \ \text{and}\ \ \det\left(\frac{\partial g}{\partial y}(x,\phi(x))\right)\neq 0\right\}
	$$
	for the function $\phi$; the authors talk of a specific solution branch $\cR$, on which they concentrate their analysis.
	Define $J(x,y)$ as the Jacobian of $(x,y)\mapsto \begin{pmatrix} f(x,y) \\ g(x,y) \end{pmatrix}$
	and define the matrices $A,B,C,D$ as functions of $x$ on the solution branch $\cR$ through the formulas
	
	\begin{equation}
		\label{ABCDdef}
		\renewcommand{\arraystretch}{1.5}
		J(x,y)=\begin{pmatrix}
			\frac{\partial f}{\partial x}(x,y) & \frac{\partial f}{\partial y}(x,y) \\
			\frac{\partial g}{\partial x}(x,y) & \frac{\partial g}{\partial y}(x,y) \\
		\end{pmatrix}
		\ \ \ \text{and}\ \ \
		\renewcommand{\arraystretch}{1}
		\begin{pmatrix}
			A & B \\
			C & D
		\end{pmatrix} := J(x,\phi(x)).
	\end{equation}
	By the chain rule,
	$$
	0=\frac{d }{d x}g(x,\phi(x)) = \frac{\partial g}{\partial x}(x,\phi(x))+ \frac{\partial g}{\partial y}(x,\phi(x))\frac{\partial \phi}{\partial x}(x)=C+D\frac{\partial \phi}{\partial x}(x),
	$$
	and hence %, additionally using the notation from \eqref{ABCDdef},
	$$
	\frac{d}{d x}f(x,\phi(x))=\frac{\partial f}{\partial x}(x,\phi(x))+\frac{\partial f}{\partial y}(x,\phi(x))\frac{\partial \phi}{\partial x}(x)=A-BD^{-1}C.
	$$
	The \emph{generalized reduced Jacobian matrix} $F_r$ with respect to the coordinate transformation $\Theta(x)$ is then defined as
	\begin{equation*}
		F_r(x,\phi(x)) := \left( \dot{\Theta}(x) +\Theta(x)  [A-BD^{-1}C]\right)\Theta^{-1}(x)%
	\end{equation*}
	with the matrix $\dot{\Theta}(x) = (\dot{\Theta}_{ij}(x))_{ij} = (\nabla\Theta_{ij}(x)\cdot f(x,\phi(x)))_{ij}$; see Section \ref{sec:diffformcontr} for similar considerations in the context of ODEs.
	Since the system \eqref{DAEgeneral} is given by $\dot{x}=f(x,\phi(x))$, on $\cR$ we have by \eqref{contraction-log} exponential contraction in $\delta z(t)=\Theta(x(t))\delta x(t)$  if
	$\mu(F_r(x(t),\phi(x(t))))\le -c<0$. For a coordinate transformation  $\rho$, the exponential contraction of
	$\delta w(t):= \rho(\phi(x(t))) \delta y(t)$, where $\delta y$ is defined through the `differential' $C\delta x+D\delta y=0$ of $g(x,y)=0$,  follows from $\sup\|\rho D^{-1}C\Theta^{-1}\|<\infty$.

	Now assume that $\cI\subset \R^n$ is a positively invariant set for system \eqref{DAEgeneral} and that the contraction condition $\mu(F_r(x,\phi(x)))\le -c<0$ holds for all $x\in \cC$, $\cI\subset \cC$; $\cC$ is called
	\emph{contraction region}. By introducing for matrices $Q\in\R^{m\times m}$ and $R\in\R^{m \times n}$  the \emph{generalized unreduced Jacobian matrix}
	\begin{equation}
		\label{unredF}
		F=
		\begin{pmatrix}
			F_r+\Theta R^TC\Theta^{-1} & \Theta R^TD\rho^{-1} \\
			Q^TC\Theta^{-1} & Q^TD\rho^{-1}
		\end{pmatrix},
	\end{equation}
	the authors show $\mu(F_r)\le \mu(F)$ for a matrix measure $\mu$ with respect to  any $p$-norm $\|\cdot\|_p$, $1\le p\le \infty$.  This implies that $\mu(F)<0$ is a sufficient condition for contraction. Note that
	$$F
	\begin{pmatrix}
		\delta z\\
		\delta w
	\end{pmatrix}= \begin{pmatrix}
		F_r\Theta \delta x +\Theta R^T(C\delta x+D\delta y)\\
		Q^T(C\delta x+D\delta y)
	\end{pmatrix}.
	$$

	To adapt this ansatz to a tractable numerical method, the authors make some simplifying assumptions.  We work the computations out in detail because there seems to be a slight error in the formulas in \cite{NVST2021contrDAE}: first, it is assumed that $\Theta\in\R^{n\times n}$ and $\rho=I\in\R^{m\times m}$ are constant, from which $F_r=\Theta[A-BD^{-1}C]\Theta^{-1}$ follows. Defining the metric $M:=\Theta^T\Theta$ and setting $R=\tilde{R}+D^{-T}B^T$ in
	\eqref{unredF}, delivers
	\begin{align*}
		F&=
		\begin{pmatrix}
			\Theta(A+ \tilde{R}^TC)\Theta^{-1} & \Theta (\tilde{R}^T+BD^{-1})D \\
			Q^TC\Theta^{-1} & Q^TD
		\end{pmatrix}\\
	&	= \begin{pmatrix}
			\Theta(A+\tilde{R}^TC)\Theta^{-1} & \Theta (B+\tilde{R}^TD) \\
			Q^TC\Theta^{-1} & Q^TD
		\end{pmatrix}
	\end{align*}
	and then
	\begin{align}
		\begin{pmatrix}
			\delta z \\
			\delta w \\
		\end{pmatrix}^T
		F\begin{pmatrix}
			\delta z \\
			\delta w \\
		\end{pmatrix}&=
		\begin{pmatrix}\delta x^T\Theta^T  & \delta y^T \end{pmatrix}
		\begin{pmatrix}
			\Theta(A+\tilde{R}^TC)\Theta^{-1} & \Theta (B+\tilde{R}^TD) \\
			Q^TC\Theta^{-1} & Q^TD
		\end{pmatrix}
		\begin{pmatrix}\Theta \delta x \\ \delta y \end{pmatrix}\n\\
		&=
		\begin{pmatrix}\delta x^T & \delta y^T \end{pmatrix}
		\begin{pmatrix}
			\Theta^T\Theta(A+\tilde{R}^TC) & \Theta^T\Theta (B+\tilde{R}^TD) \\
			Q^TC & Q^TD
		\end{pmatrix}
		\begin{pmatrix} \delta x \\ \delta y \end{pmatrix}\n\\
		&=
		\begin{pmatrix} \delta x \\ \delta y \end{pmatrix}^T
		\underbrace{\begin{pmatrix}
				M & 0 \\
				0 & I
			\end{pmatrix}
			\begin{pmatrix}
				I & \tilde{R}^T \\
				0 & Q^T
		\end{pmatrix}}_{=:Z^T}
		\underbrace{
			\begin{pmatrix}
				A & B \\
				C& D
		\end{pmatrix}}_{=J(x,\phi(x))}
		\begin{pmatrix} \delta x \\ \delta y \end{pmatrix}\n.
	\end{align}
	Note that the matrix $Z$ is affine in the components of the metric $M=\Theta^T\Theta$, and  one can get an appropriate $\Theta$ from $M$ by the Cholesky decomposition.
	With $\alpha=(x,y(x))$, the authors then consider the system \eqref{DAEgeneral} close to an equilibrium $\alpha^*$ and assume that  $J(\alpha)=J^*+\sum_k\alpha_k J_k$ is a reasonable approximation, i.e.~that
	$J(x,\phi(x))=J(\alpha)$ is affine in the components $\alpha_k$ corresponding to the components of $x$ and $y$. Then $\mu(F)\le -\beta<0$, where $\mu$ is the matrix measure with respect to the Euclidean norm $\|\cdot\|_2$,
	can be written as the bilinear matrix inequality (BMI)
	\begin{equation}
		Z^TJ(\alpha)+J^T(\alpha)Z  \le -\beta I.
	\end{equation}
	The paper is concluded by discussing how to generate positively invariant sets around $\alpha^*$ and use numerical methods for BMIs, e.g.~from \cite{KZM2018BMI1,KZM2018BMI2}, to generate contraction metrics for concrete systems, which is then used to show contraction of a two-bus system.

	\subsubsection{Continuous, piecewise affine contraction metrics}

	Given a \emph{shape regular triangulation} $\cT=\{\fS_\nu\}$ of a compact set $\cD\subset \R^n$, i.e.~a subdivision of $\cD$ into simplices $\fS_\nu$ that pairwise intersect in a common face or not at all, one can parameterize
	\emph{continuous, piecewise affine} (CPA) functions by fixing their values at the vertices $\cV$ of the simplices.  Alternatively, one can approximate a function $g\colon \cD \to \R$ by interpolating its values over the simplices,
	using its values at the vertices of the simplices; called \emph{CPA interpolation} and denoted by $g_{\operatorname{CPA}}$.
	This has been used to compute CPA Lyapunov functions for nonlinear systems using linear programming, see e.g.~\cite{JGD1999cpajulian,Mar2002cpa,GiHa2014CPArev}, and has more recently
	been adapted to compute CPA contraction metrics using LMIs.  The idea is to construct a contraction metric $M\colon \cD \to \cS_n$
	that is a CPA function. To this end, inequalities for the conditions of  a contraction metric in the values $M_x\in \cS$ at the vertices $x\in \cV$ are derived.
	%, taking the approximation errors in $g_{\tiny\operatorname{CPA}}$ into account.
	This results in  LMIs and one can incorporate
	estimates for each simplex % approximation errors
	in
	such a way that $M$ parameterized by a feasible solution
	%, that delivers the values $M_x$,
	is a true contraction metric for the system.

	Let us explain this in more detail.  Each simplex $\fS_\nu=\conv\{x^\nu_0,\ldots,x^\nu_n\}$ is the convex hull of $(n+1)$ vertices in $\cV$ and every $y\in \fS_\nu$ can be written uniquely as a convex combination
	%$$
	%y=\sum_{k=0}^n\lambda_k x^\nu_j,\ \ \lambda_k\ge 0,\ \  \sum_{k=0}^n\lambda_k =1
	%$$
	of the vertices. If we assign to each vertex $x\in \cV$ a value $M_x=([M_{ij}]_x)\in\cS_n$, we can define a function $M\colon \cD\to\cS_n$ through
	$$
	M(y)=\sum_{k=0}^n\lambda_k M_{x^\nu_k},\ \ \ \text{where}\ \ \ y=\sum_{k=0}^n\lambda_k x^\nu_k,\ \ \lambda_k\ge 0,\ \  \sum_{k=0}^n\lambda_k =1.
	$$
	Each entry $M_{ij}$ of the function $M=(M_{ij})$ is continuous and affine on each simplex in $\cT$.  In particular, $M_{ij}$ can be written as $M_{ij}(y)=\left(\nabla [M_{ij}]_\nu\right)^T( y-x^\nu_0) + [M_{ij}]_{x^\nu_0}$ on $\fS_\nu$, where the vector $\nabla [M_{ij}]_\nu\in\R^n$ is constant on the simplex $\fS_\nu$ and  linear in the values $[M_{ij}]_{x^\nu_k}$.
	%For the CPA interpolation $g_{\small\operatorname{CPA}}$ of the function $g$ we can bound the approximation error by
	%\begin{equation}
	%\label{cpaappr}
	%|g(y)-g_{\operatorname{CPA}}(y)| \le h^2 B_g,
	%\end{equation}
	%where $h$ is an upper bound on the diameters of the simplices in $\cT$ (in the Euclidian norm) and $B_g$ is an appropriate upper bound on the second-order derivatives of $g$ on $\cD$.
	%

	In Giesl \& Hafstein 2013 \cite{GiHa2013CPAmetric}, this was used to derive LMIs, or a semidefinite optimization problem, for systems $\dot{x}=f(t,x)$, where $f$ is $T$-periodic, i.e.~$f(t+T,x)=f(t,x)$ for all $(t,x)\in\R\times \R^n$.  A feasible solution of the  semidefinite optimization problem can be used to parameterize a contraction metric for the system that asserts the existence and exponential stability of a periodic solution, which is a periodic orbit on the manifold $S^1_T\times \R^n$, where $S^1_T$ is the circle of circumference $T$. Further, it is shown that if the system has an exponentially stable periodic solution, then the optimization
	problem always has a feasible solution if the domain under consideration is subdivided into sufficiently small  simplices.
	
	In Hafstein \& Kawan 2019 \cite{HK2019}, a similar approach is followed, but additionally combining methods from the dimension estimation of attractors, to obtain upper bounds on the \emph{restoration entropy} of nonlinear systems.
	Note that semidefinite optimization is not as mature a  subject as linear programming and of higher computational complexity.  In particular, large semidefinite optimization problems cannot be solved efficiently with state-of-the-art methods and hardware.
	Therefore, a simplified procedure is proposed, computing a constant metric in the first optimization
	problem.  Then the second optimization problem simplifies to a linear programming problem and this method is used to obtain an upper bound on the restoration entropy for the Lorenz system.
	{In more detail, first} a metric is computed, which is subsequently fed into a different semidefinite optimization problem that delivers a Lyapunov-like function, and together the metric and the Lyapunov-like function deliver an upper bound on the restoration entropy.

	In Giesl, Hafstein \& Mehrabinezhad 2021 \cite{iman-eq}, a different strategy is followed to compute contraction metrics for nonlinear systems with an exponentially stable equilibrium.  First, a semidefinite optimization problem is proposed, of which every feasible solution delivers a contraction metric for the system.  The essential constraints of the problem for a system $\dot{x}=f(x)$ are of the form
	\begin{equation}
		\label{cparbfcollo}
		\frac{\partial f^T}{\partial x}(x^\nu_j)M_{x^\nu_j}+ M_{x^\nu_j}\frac{\partial f}{\partial x}(x^\nu_j)+\left(\nabla [M_{ij}]_\nu\cdot f(x^\nu_j)\right)+h_\nu^2E_\nu I \le -\varepsilon I
	\end{equation}
	for all vertices $x^\nu_j$, $j=0,\ldots,n$, of all simplices $\fS_\nu$ in the triangulation.
	Here, $h_\nu$ is the diameter of the simplex $\fS_\nu$ and $E_\nu$ is an error term involving an upper bound on $M$ on $\fS_\nu$ and $\|\nabla [M_{ij}]_\nu\|_1$, which are implemented through auxiliary variables.
	Instead of solving this optimization problem, collocation as in Section \ref{sec:collo} is used to compute values for the variables $M_x$ and then the constraints \eqref{cparbfcollo} are verified.
	If all constraints are fulfilled, then the CPA interpolation of the metric computed by collocation is a contraction metric for the system.
	This results in a method that combines the numerical efficiency of the collocation method (solving equations) with the rigour of the LMI approach.  Further, it is proved that this method is able to compute a contraction metric for any system with an exponentially stable equilibrium,
	given that one uses a sufficiently dense collocation grid and sufficiently small simplices.
	In \cite{iman-per}, this method is adapted to compute contraction metrics for systems with exponentially stable periodic orbits and in \cite{GHM2021rob} to demonstrate the robustness of the metrics computed by this method.

	%Verification of computed contraction metric via RBF for equilibrium \cite{iman-eq} and for periodic orbit \cite{iman-per}

	\subsection{Subgradient method}
	
	A subgradient algorithm can be used to estimate the restoration entropy, see Kawan, Hafstein \& Giesl 2021 \cite{KHG2021}. According to \cite{kawan2021remote}, the restoration entropy of a system given by the autonomous ODE $\dot{x} = f(x)$ with $f \in C^1(\R^n,\R^n)$, on a compact positively invariant set $K \subset \R^n$, satisfies
	\begin{equation}\label{eq_resent}
		h_{\mathrm{res}}(f,K) = \inf_{M \in C^1(K,\mathcal{S}^+_n)} \max_{x\in K} \max_{0 \leq k \leq n} \sum_{i=1}^k \zeta_i^M(x)
	\end{equation}
	with $\zeta_1^M(x) \geq \zeta_2^M(x) \geq \cdots \geq \zeta_n^M(x)$ being the solutions $\zeta(x)$ of
	\begin{equation}\label{eq_ct_alg_eq}
		\det\left[
		\frac{\partial f^T}{\partial x}(x)M(x)+ M(x)\frac{\partial f}{\partial x}(x)+\dot{M}(x) - \zeta (x)M(x) \right] = 0.
	\end{equation}
	In the numerical approach to solve the minimization problem posed by \eqref{eq_resent}, only Riemannian metrics on $K$ of the form $$M(x)=e^{p_a(x)}P$$ are considered, where $P\in{\mathcal S}_n^+$, $p_a(x)=\sum_{\alpha\in I}a_\alpha x^\alpha$ is a polynomial and $I\subset \mathbb N_0^n$ a finite set of multi-indices. Hence, the search for the optimal metric is restricted to the space $\R^{|I|}\times {\mathcal S}^+_n$, where $\mathbb R^{|I|}$ is equipped with the usual Euclidean metric and ${\mathcal S}^+_n$ with its standard Riemannian metric, the so-called trace metric given by
	\begin{equation*}
		\langle V,W \rangle_P := \tr(P^{-1}VP^{-1}W) \mbox{\quad for all\ } P \in {\mathcal S}_n^+,\ V,W \in T_p{\mathcal S}_n^+ = {\mathcal S}_n.
	\end{equation*}
	The optimization is performed using the Riemannian subgradient algorithm for geodesically convex functions on the Riemannian product manifold $\cM:=\R^{|I|} \times {\mathcal S}^+_n$, cf.~\cite{FOl,Fea}. The required geodesic convexity of the objective function
	\begin{equation*}
		(a,P) \mapsto \max_{x\in K} \max_{0 \leq k \leq n}\sum_{i=1}^k\zeta_i^{e^{p_a(\cdot)}P}(x)%
	\end{equation*}
	for this approach is shown in \cite{KHG2021}.
	
	For the algorithm, one first fixes an initial metric, e.g.~$M_0(x)=I$ (i.e.~$a_0 = 0$ and $P_0=I$), and a sequence of step sizes $t_j$ satisfying $\sum_{j=1}^\infty t_j = \infty \mbox{\quad and \quad} \sum_{j=1}^\infty t_j^2 < \infty$. Then, in the $j$-th step of the algorithm, one computes a subgradient $S_j$ of the objective function at the current point $(a_j,P_j)$, normalizes $S_j$ to unit length by putting $\bar{S}_j := S_j/|S_j|$, and then updates the point (i.e.~the Riemannian metric) by
	\begin{equation*}
		(a_{j+1},P_{j+1}) := (a_j - t_j \bar{S}_j^1,\exp_{P_j}(-t_j \bar{S}_j^2)).
	\end{equation*}
	Here, $\exp_{P_j}(\cdot)$ is the Riemannian exponential map of ${\mathcal S}^+_n$ at $P_j \in {\mathcal S}^+_n$, and $\bar{S}_j^1 \in \R^{|I|}$, $\bar{S}_j^2 \in T_{P_j}{\mathcal S}^+_n = {\mathcal S}_n$ are the two components of $\bar{S}$. To compute the required subgradient $S_j$, one first determines a point $x^*\in K$ that maximizes
	\begin{equation}\label{eq_ct_maximization}
		g(x;M_j) := \max_{0 \leq k \leq n} \sum_{i=1}^k \zeta_i^{M_j(\cdot)}(x),
	\end{equation}
	{where $M_j(x)=e^{p_j(x)P_j}$,}
	and then computes a subgradient of the function $(a,P) \mapsto \max_{0 \leq k \leq n} \sum_{i=1}^k \zeta_i^{e^{p_a(\cdot)}P}(x^*)$ at $M_j$. The computed subgradient is then also a subgradient of the objective function. Each $b_j := g(x^*;M_j)$ is an upper bound on the restoration entropy. Note that, in general, the subgradient algorithm does not deliver a monotonically decreasing sequence $b_j$, but it can be shown that $c_j:=\min_{0\le i \le j} b_i$ converges to the infimum of the objective function restricted to $\cM$ and the corresponding $M_j$ converge to a minimizing metric, if such a minimizer exists. If no minimizer exists, then still $\lim_jc_j=\liminf_j b_j$ is equal to the infimum of the restricted objective function.
	
	This method was used in \cite{KHG2021} to compute estimates of the restoration entropy for continuous- and discrete-time systems, and formulas for the subgradients as discussed above are provided. A description of the software can be found in Kawan, Hafstein \& Giesl 2021 \cite{KHG2021X}.
	
	In a similar way, the subgradient algorithm can be used for the computation of contraction metrics for equilibria and periodic orbits as well as the estimation of the dimension of attractors. In the case of the equilibrium {\cite{Magnea}}, e.g., \eqref{eq_ct_maximization} is replaced by the maximal {$\zeta^{M_j(\cdot)}$}, namely
	\begin{equation}\label{eq-subgradient}
		g(x;M_{j}) := \zeta_{1}^{M_j(\cdot)}(x).
	\end{equation}
	
	{Due to the maximization involved in the computation of a subgradient, the subgradient method in the form described here suffers drastically from the curse of dimensionality. However, for higher dimensional systems it is very likely that a small-gain approach as developed in \cite{matveev2019comprehending} can be used to reduce the complexity by exploiting a decomposition of the system into low-dimensional interconnected subsystems and performing the costly computations only for those low-dimensional systems.}
	
	%
	%
	%Further literature
	%\begin{itemize}
	%	\item B. P. Demidovich, "Dissipativity of a nonlinear system of differential equations, Part II," Vestnik Moscow State Univ. Ser. Matem. Mekh., vol. 1, pp. 3–8, 1962.
	%	\item Krasovskii Stability of Motion 1963 book in library
	%\end{itemize}

	\section*{Concluding Remarks}
	
	In this review, we have presented the basic ideas behind contraction analysis, its connection to the Lyapunov stability theory, and its historical development.  Further, we have discussed the different criteria for contraction and convergence, the relationships between them, and various extensions of the theory to different classes of systems. Finally, we have described how collocation, semidefinite optimization, and the subgradient method on matrix manifolds, can be applied to compute contraction metrics for various kinds of systems. The computational methods are of high complexity and thus currently only capable of computing contraction metrics in low-dimensional systems. However, this could be improved in the future by either using small-gain approaches, improved numerical methods or other developments.
	It is the hope of the authors that this review is useful for researchers interested in contraction analysis, its applications to real-world systems, and numerical methods for their computation.

%The acknowledgments section should not be numbered. Please begin it with \section*{Acknowledgments}
%\section*{Acknowledgments} We would like to thank you for \textbf{following
%the instructions above} very closely in advance. It will definitely
%save us lot of time and expedite the process of your paper's
%publication.

% You may incorporate your references as follows in your main tex file.
% Using BibTex is not recommended but can be handled. If you use BibTex, please include the file with your final paper files.
% AIMS editorial staff will add MR and DOI numbers to your references.

{\small
	\bibliographystyle{plain}
	\bibliography{LibContMetr-new2}
}

\medskip
% The information below will be filled in by AIMS editorial staff
Received 8 October 2021; revised 15 February 2022; early access xxxx 20xx.
\medskip

\end{document}